\documentclass [12pt, book, reqno, a4paper, czech, french, hyperref]{amsart}
\usepackage[dvips]{epsfig}
\usepackage{graphicx}
\usepackage[dvips]{graphics}

\usepackage{times}

\usepackage{color}
\def\bl{\color{black}}
\def\red{\color{black}}

\newcommand{\cqfd} {\hbox {\unskip \kern 6pt \penalty 500\raise -2pt \hbox
{\vrule \vbox to 5pt {\hrule width 4pt \vfill \hrule }\vrule
}\par }}

\newlength\jataille 
\newcommand{\figgauche}[3]%
{
\jataille=\textwidth\advance\jataille by -#1
\advance\jataille by -.5cm
\begin{minipage}[a]{#1}
\includegraphics[width=#1]{#2}
\end{minipage}
\vskip2mm
\begin{minipage}[a]{\jataille}
\footnotesize #3 \normalsize
\end{minipage}
}

\usepackage{amsfonts}
\font\tencyr=wncyr10          
\def\cyr{\tencyr
}

\def \Ee{\char"03 }
\def \mz{\char"7E }

\def\ppref
#1 
{$\bullet $
\ref{#1}
}
\newcommand{\lettre}[1]{
\refstepcounter{section}
\vskip 5mm
\centerline{Lettre \thesection}
\sectionmark{Lettre \thesection}
\vskip 3mm }

\usepackage{ulem}
\usepackage[francais]{babel}

\usepackage{amsfonts}
\font\tencyr=wncyr10          
\def\cyr{\tencyr
}

\def \Ee{\char"03 }
\def \mz{\char"7E }

\def \I'{\accent"24 I}

\begin{document}

\Large 
\centerline{\red \bf The emergence of French probabilistic statistics}
\centerline{  \red \normalsize \it Borel and the Institut Henri Poincar\'e around the 1920s}
\vskip 3mm 
\large 
\centerline{R\'emi CATELLIER\footnote{ENS Cachan, {France}. remi.catellier@gmail.com }
and 
Laurent MAZLIAK\footnote{Laboratoire de Probabilit\'es et Mod\`eles al\'eatoires, {Universit\'e Pierre et Marie Curie, Paris, 
France}. laurent.mazliak@upmc.fr . The paper was completed when this author participated to a research program supported by the Institute of Communication Sciences of the CNRS.}
}
\normalsize 
\centerline{\it 
 }
\renewcommand{\rightmark}{\uppercase{\bf The emergence of French probabilistic statistics }}
\renewcommand{\leftmark}{{\bf{R\'emi Catellier  and Laurent Mazliak}}}
\vskip 5mm 

\centerline{\bf Abstract}

\begin{quote}{\tiny This paper concerns the emergence of modern mathematical statistics in France after  the First World War. Emile Borel's achievements are presented, and especially his creation of two institutions where mathematical statistics was developed: the {\it Statistical Institute of Paris University}, (ISUP) in 1922 and above all the {\it Henri Poincar\'e Institute} (IHP) in 1928. At the IHP, a new journal {\it  Annales de l'Institut Henri Poincar\'e}  was created in 1931. We discuss the first papers in that journal dealing with mathematical statistics.}
\end{quote}

\section*{Introduction}

The important transformations in the field of mathematics of randomness between 1900 and 1930 are now rather well understood. Several large-scale studies have been published which present wide pictures of the actors and ideas involved in what may be considered as a major evolution in the scientific life of the 20th Century. See in particular several chapters in the two volumes of the the monumental treatise \cite{Krugerandal1987} and fundamental book \cite{VonPlato1994}. These books have been complemented by many recent papers concerning more detailed aspects of this story (see among numerous others \cite{Bru2003}, \cite{HavlovaMazliakSisma2005}, \cite{SiegmundSchulze2006}
 and other references in these papers). Yet it seems that these studies are often more specifically centered on the  probabilistic aspects of the question than on the statistical side. When one consults the recent, very comprehensive collective work on statisticians edited by Heyde and Seneta \cite{HeydeSeneta2003}, it is striking to see that many of those who are mostly known today as specialists in probability theory were also involved in the shaping of modern mathematical statistics. This includes for example the Austro-German Richard von Mises, the Soviets E.B. Slutzky, A.Y.Khin\v cin and A.N. Kolmogorov,  the Italian Guido Cantelli and naturally, the French Emile Borel about whom we shall write at length below. 
Though basic tools of the modern theory of probability (especially the use of Borel's measure theory of sets and the Lebesgue integral) had been available since the 1910s in France, it took a long time for  mathematics of randomness, above all the most modern ones, to penetrate the quite reluctant and suspicious world of field statisticians who preferred the use of descriptive methods without deep mathematical theory. The present paper approaches how eventually mathematical statistics slowly emerged in the 1920s. Let us mention that this question has already been studied, especially by Desrosi\`eres in his impressive book \cite{Desrosieres2000} where a 
wide panorama of the international situation is provided. More particularly, in France, we can also mention among others the paper \cite{Meusnier2006}. Our paper seeks to provide further information on the  French situation and its relations with what was happening abroad, and how the link was finally drawn between probability and statistics. 

As already mentioned, Emile Borel (1871-1956) was a central actor in this transformation, as were his successors Maurice Fr\'echet (1878-1973) and  Georges Darmois (1888 -1961). By means of original pedagogical and scientific initiatives, the three men participated in the setting up  of two institutes in Paris, the {\it Institut de Statistiques de l'Universit\'e de Paris} (ISUP) in 1922, and above all the {\it Institut Henri Poincar\'e} (IHP) in 1928 where mathematical statistics was presented and recent findings in the subject described for the first time in France. {\red A good hint of the role played by the aforementioned mathematicians is that they were, between 1920 and 1950, the only mathematicians elected as President of Paris Statistical Society (SSP) (see \cite{Mazliak2010}).}

Our work stems from the conference organized in November 2008 by Mich\`ele Audin and Catherine Goldstein for the 80th anniversary of the IHP. The second author was invited to present a talk on the same theme as the present article. Moreover, this talk was largely inspired by a short memoir \cite{Catellier2008} written by the first author after a short research stay at University Paris 6 during Spring 2008. On that occasion, important new sources on the subject (archival material and various articles) had been collected. This work was likewise strongly connected to a large international research program about the mathematics around the Great War, initiated 6 years ago. This program opened up numerous new directions for the historiography of this period, which were in particular discussed during an  international conference in 2007 \cite{Luminy2007}.

{\bl 
In 1996, Stigler (\cite{Stigler1996}) proposed the somewhat provocative hypothesis that Mathematical Statistics began in 1933. Some facts presented in our paper may therefore appear as a contradiction to Stigler's thesis. This is in fact only partly the case.  Reading Stigler's paper shows that the author concentrates mainly on the case of the United States; the obvious predominance of Anglo-Saxon statistics after the 1940s justifies particular attention on what happened in the USA during the interwar period. When we decided to scrutinize the French situation in the same period, we knew that we were looking at a country in which the role of mathematics in statistical studies was years behind in comparison with other countries. However, the actors of our story (Emile Borel and his successors) were mathematicians who used mathematics to teach and (timidly) invent statistical methods beginning already in the 1920s. As we shall see, it is clear that in their mind they were establishing a discipline called {\it mathematical statistics} in France - and {\it Mathematical Statistics} is precisely the title of Darmois' book in 1928 \cite{Darmois1928}.  In fact, the point seems mainly to be a question of scale:  Stigler's hypothesis pertains to a more macroscopic vision, so to speak. He himself mentions in his paper that it would certainly be more accurate to propose a confidence interval than to specify a precise year such as 1933. The fragments of the history of statistics we shall deal with in this paper may be considered as a part of the powerful stream that Stigler emphasizes, a stream which established mathematical statistics within the mathematical sciences. }

{\red 
The present paper is therefore an attempt to prove that Borel and his successors Darmois and Fr\'echet wanted to develop a type of statistics based on theoretical results of the calculus of probability. In their mind, the relationship between mathematical statistics and probability was therefore a {\it filiation} more than a complementarity. Moreover Fr\'echet, in his presidential address when he became President of the SSP in 1948, explicitly mentioned this idea and called calculus of probability the {\it father of mathematical statistics} (\cite{Frechet1948}). The three parts of our article discuss how the aforementioned program was implemented.  

 In the first part, we shall try to draw a picture of Emile Borel and his encounter with randomness. There are several texts presenting Borel and probability, and notably the important study by Callens \cite{Callens1997}. However, it seems necessary to recall several facts about this major actor of our story in order to explain how his scientific and social personality contributed to the developments we discuss. Of particular importance for our subject were Borel's unceasing efforts to understand when and how probability could be applied legitimately  to statistical situations.  We believe that the outlines of the mathematics taught in the new institutions during the 1920s and 1930s were strongly related to Borel's point of view. We insist on the role played by Borel's activities during the Great War. They provided him with a large-scale experience of applied statistics, especially through his various contacts with politicians. 
 
 The second part is devoted to the institutional realization of Borel's program. We deal with Borel's fundamental belief in the role of mathematical statistics in the education of decision-makers. We present how the setting up of the ISUP, and later of the IHP, was for him an opportunity to give shape to his ideas on education and scientific politics. We will dwell in particular on the roles played by Fr\'echet in Strasbourg and  Darmois in Paris in the emergence of mathematical statistics during the 1920s in France. We will also provide detail on the functioning of the IHP and its scientific program.

 The third part focusses on how the aforementioned institutions were used to implement what can be described as a {\it technology transfer} aimed at filling a gap in French statistical techniques by means of importation of foreign achievements. British, Scandinavian, German and Italian research was especially scrutinized. After having commented on the creation of the {\it Annales de l'Institut Henri Poincar\'e}, we concentrate on the first two papers on statistics published in the journal in 1932, the first one by Darmois, the second by Guldberg. The study of these papers,  with a mathematical insight into their contents, is used as a demonstration of how the technology transfer worked. Moreover, the selection of subjects and speakers at the IHP shows how  a kind of mathematical statistics {\it \`a la fran\c caise} emerged based on 
the know-how of analysts and mathematical physicists converted to probability. It expanded during the 1930s as illustrated by Fr\'echet's polemics at the International Statistical Institute about the correlation (see \cite{Armatte2001}). }

\section{Emile Borel's apprenticeship in randomness}
Emile Borel, one of our main protagonists, was a major figure of the French mathematical scene in the first half of the 20th Century. His works were published under the care of his friend and former student Maurice Fr\'echet (\cite{Borel1972}). On this occasion, Fr\'echet wrote a first biographical sketch, expressing his deep admiration for his former master. Since then, several biographies have been published about him. Let us note \cite{Politchuk1973} and more recently \cite{Guiraldenq1999}. 

Borel's mathematical achievements are so rich and offer so great a variety of subjects that it is obviously not the place here to propose a comprehensive picture of them. Besides, there is now a collection of works dealing with several aspects of Borel's mathematical production which will give the interested reader an idea not only of the results obtained by Borel but also of the spirit in which he produced them. We note for example the book \cite{Hawkins1970} in which the genesis of the theory of measurable sets is presented in depth (a theory on which Lebesgue based the construction of his integral), or the paper \cite{MaureyTacchi2005} dealing with Borel's covering theorem (the compactness of the interval $[a,b]$). The very recent paper \cite{Barberousse2009} deals with philosophical aspects of Borel's consideration on approximation and with his criticisms of a notion of idealistic knowledge (a subject deeply related to ours as we shall see below). 


The object of this first section is therefore to briefly recall how Borel became familiar with probability, and how randomness (in particular through a statistical orientation) acquired a considerable dimension in his  vision of a scientifically enlightened citizenry and the rational direction of state business. 

\subsection{Biographical sketch}

Emile Borel was born to a Protestant middle-class  family in Saint-Affrique, in Aveyron,  in the center of southwest France.  He kept close ties with Saint-Affrique throughout his life.  After brilliant secondary studies, he went to Paris to sit for the competitive examinations leading to the Grandes Ecoles, the schools where the French scientific and administrative elites are trained.  There he studied under the famous teacher Boleslas Niewenglowski, along with the son of mathematician Gaston Darboux, and he later recounted that it was at Darboux's home that he discovered his passion for scientific and especially mathematical research.  The Ecole Normale Sup\'erieure was the place where he pursued this passion.

However, purely mathematical achievements were only one side to the rich personality of Emile Borel. Very young, he engaged in an original mixture of academic, political and philosophical disputes, and had the strong desire not to be one of the type of mathematicians locked in university libraries, speaking only to his peers in his ivory tower. Borel was always convinced that the scientist had a role to play in the society of his time alongside politicians, writers, and artists, and he had the strong desire to  participate in such social exchanges. Borel engaged himself on the social stage, especially through circles connected to the family of his talented wife Marguerite and his father-in-law, mathematician Paul Appell.  Marguerite wrote fiction under the pen name Camille Marbo (for MARguerite BOrel); in 1913, she won the Femina prize for her novel {\it La statue voil\'ee}.  In 1905, Borel and Marbo founded a monthly journal, the {\it Revue du Mois}, which for 10 years was a leading general intellectual outlet for the moderate French left. When he became a politician in the 1920s, Borel, with his Protestant radical-socialist background, was passionately attached to secularism and rationalism,   and was a typical representative of the {\it republic of professors} as journalist Albert Thibaudet later described the 3rd Republic\cite{Thibaudet1927}. However, in 1905, his ambition was at first only intellectual and the {\it Revue du Mois} became his favorite place for presenting his ideas (and his friends' views) on every kind of subject. Let us recall besides that this aspect of Borel's personality was not always accepted by his colleagues. Lebesgue, in particular, who had no great taste for the social commitments of scientists, criticized his behavior. 

\begin{quote}
To speak the truth I reproach you for the Revue du Mois. I know quite well that you find there an opportunity  to put into practice your qualities of action and your zeal as an administrator, but it is what I esteem least in you\footnote{Pour tout dire je vous reproche La Revue du Mois. Je sais bien que vous trouvez l\`a l'occasion de d\'epenser vos qualit\'es d'action et vos ardeurs d'administrateur, mais c'est ce que j'estime le moins chez vous.}. 
({\it Lebesgue to  Borel, [20] February 1909} \cite{Dugac2005})
\end{quote}
The interested reader can find details about the foundation of the journal in \cite{Mazliak2007}.

For several reasons, mixing an original turn in his mathematical reflection with his thoughts about how to enlighten his fellow man on a scientific basis, Borel soon had a kind of revelation that the mathematics of randomness could be an essential field for this struggle. Contrary to what is still sometimes asserted, Borel's reflections on probability do not belong to a second (and much less brilliant) scientific life  beginning roughly in the 1920s. Borel met probability at the very beginning of the century and the subject always remained on his mind. He developed a complex relationship  to it, one not always perfectly understood by his colleagues and in particular by younger French mathematicians, such as Paul L\'evy,  who became involved in the probabilistic field. Anyway, there is general agreement today that Borel should be seen as  the harbinger of a renewal mathematical probability at the beginning of 20th century, opening the path to the axiomatic formalization based on measure theory achieved by Kolmogorov in his {\it Grundbegriffe der Warscheinlichkeitsrechnung} \cite{ShaferVovk2006}.

Several texts deal with Borel's original attitude towards randomness. Let us quote in particular \cite{Knobloch1987} and obviously the seminal book \cite{VonPlato1994} where Borel's probabilistic turn is studied in depth. The paper \cite{BruBruChung1999} provides another picture of this turn. Finally, in the second part of his PhD (published as \cite{Callens1997}), Callens offers a profound examination of Borel's considerations on the role of quantified randomness in various aspect of social life. Much of what we write about Borel's life is inspired by these texts. 

\subsection{Randomness}

When Borel entered the Ecole Normale, he immediately specialized in mathematics, beginning fundamental studies on divergent series, for which he introduced different modes of summability.  This soon led him to fundamental work on the measure of sets, which cleared the way for Lebesgue to construct his integral and revolutionize analysis \cite{Hawkins1970}.  Measure theory also led Borel to focus on probability theory, starting in 1905 in \cite{Borel1905} . Though we do not know exactly how Borel realized that measure theory and Lebesgue integral were perfectly adapted to probabilistic considerations, we can reasonably find two converging sources for his interest for this field. A first one is Wiman's paper \cite{Wiman1900}  where, as Borel himself mentions in \cite{Borel1905}, measure theory was applied for the first time in a probabilistic context to obtain a limit distribution of the quotients of the continued fraction representation of a real number chosen at random between 0 and 1 (see details in \cite{Mazliak2009}).  It is worth recalling that Borel had introduced considerations on the measure of the sets in the first place (in his studies of analytic continuation \cite{Borel1894} - see \cite{Hawkins1970} pp. for details) as a powerful tool for proving the existence of a mathematical object. Borel's first reasoning in that direction was to assert that the complement of a negligible subset (i.e. one with null measure) of an interval must contain at least one point. As Boutroux wrote \cite{Boutroux1920}, Borel was a  figure departing from a constructivist ideal, and he probably had the striking intuition that a non-constructivist attitude was extremely well adapted to describe mathematically the vague notion of randomness. In his thesis {\it Essay on Approximate Knowledge} (Essai sur la connaissance approch\'ee), published in 1927 \cite{Bachelard1927}, Bachelard saw the use of probability (in particular by Borel) as an adequate tool for abstract speculations because probability is {\it a notion refractory to any ontology}. However, he was not completely right concerning Borel, for whom a  reasoning involving probability provided a very concrete condition of existence, based on the precise meaning Borel had given to the notion of probability. As Darmois later wrote in  \cite{Darmois1958}, {\it until the beginning of XXth century, the very notion of probability remained rather vague. It is to Emile Borel that we owe its rigorous definition, based on the measure of sets\footnote{Jusqu'au d\'ebut du XX\`eme si\`ecle, la notion m\^eme de probabilit\'e restait assez vague. C'est \`a Emile Borel que l'on doit sa d\'efinition rigoureuse, bas\'ee sur la mesure d'ensemble}.}

But obviously the most convincing reason for Borel to be interested in probability should have been the increasing presence of probabilistic considerations in the new physics of his time, above all in statistical mechanics. Borel appears in this regard as a successor of Poincar\'e, who, after having tried to circumvent the cumbersome and shocking presence of randomness in physics, had finally made the best of the bad situation and decided, as he was no more in a position to avoid it, to give it a reasonable mathematical shape. After 1890, considerations on randomness and probability were always present in Poincar\'e's mind and he wrote the textbook \cite{Poincare1896} for his students at the Sorbonne, as well as several texts aimed at a more general audience published in his books (see in particular Science and Hypothesis \cite{Poincare1902} ). 

However, Borel (who belonged to the next generation) adopted a slightly different viewpoint from Poincar\'e's . For the latter  what was at stake was to `convince himself'   that  the probabilistic hypothesis was  tolerable in physics (as Bru wrote in \cite{Bru2003}). 

As Barberousse comments in her paper \cite{Barberousse2009}, it seems that for Poincar\'e probabilsm was reasonable only if it had a provisory character. For Borel, 
on the contrary, some situations (as the one he considers in his paper on the paradox of the weat heap - see section \ref{polemics} below) generate questions for which the adapted answer can only be given using a probabilistic modelization. It is therefore interesting to develop the mathematics of randomness in order to provide new tools for elaborating and 
studying these models.

Borel's first important probabilistic paper \cite{Borel1906}, devoted to a mathematical presentation for Maxwell-Boltzmann's kinetic theory of gases, was a magisterial illustration of such a program. He wrote

\begin{quote} I would like to address all those who shared Bertrand's opinion about the kinetic theory of gases, that the problems of probability are similar to the problem of finding the captain's age when you know the height of the mainmast. If their scruples are partly justified because you cannot blame a  mathematician with his love of rigor, it nevertheless does not seem to me impossible to make them happy.
This is the aim of the following pages : they do not bring any real advance in the theory from the physical point of view; but perhaps they will result in convincing several mathematicians of its interest, and, by increasing the number of researchers, will indirectly contribute to its development. If this is the case, they will not have been useless, independently of the esthetic interest connected with any logical construction\footnote{Je voudrais m'adresser \`a tous ceux qui, au sujet de la th\'eorie cin\'etique des gaz, partagent l'opinion de Bertrand que les probl\`emes de probabilit\'e sont semblables au probl\`eme de trouver l'\^age du capitaine quand on conna\^\i t la hauteur du grand m\^at. Si leurs scrupules sont justifi\'es jusqu'\`a un certain point parce qu'on ne peut reprocher \`a un math\'ematicien son amour de la rigueur, il ne me semble cependant pas impossible de les contenter.
	C'est le but des pages qui suivent : elles ne font faire aucun progr\`es r\'eel \`a la th\'eorie du point de vue physique; mais elles arriveront peut \^etre \`a convaincre plusieurs math\'ematiciens de son int\'er\^et, et, en augmentant le nombre de chercheurs, contribueront indirectement \`a son d\'eveloppement. Si c'est le cas, elles n'auront pas \'et\'e inutiles, ind\'ependamment de l'int\'er\^et esth\'etique pr\'esent dans toute construction logique}.\end{quote} 

Borel's most important probabilistic achievement was the seminal paper \cite{Borel1909} devoted to denumerable probabilities and how they could be used for arithmetical considerations. It remains probably one of the most important probabilistic papers ever done. 
The study of continued fractions, which had considerable consequences for modern probability theory - see \cite{VonPlato1994}, \cite{Mazliak2009}, appeared there once again. Borel stated there a first form of the strong law of large numbers, catching his contemporaries {\it by surprise} (\cite{VonPlato1994}, p.57) with this strange way of proving that almost every real number satisfied a property though it was hard to decide whether any particular number satisfied it. In fact, for Borel, such an assertion was a brilliant justification for his aforementioned intuition that a probabilistic answer was well adapted for providing a new modality of existence. The central point was to describe an impossibility by asserting that a set is negligible. This situation of almost impossibility or almost certitude was the only one in which one may attribute an objective value to probability. For Borel, the most important notion of probability theory (he would even call it later {\it La loi unique du hasard}, the unique law of randomness) was the interpretation that events with minute probability were impossible, as was illustrated by his famous popular image of the typing monkey reconstituting the books of the National Library. For Borel, the calculus of probability was seen as an application of mathematical analysis. One must therefore consider its results with the same caution as for any other application of mathematics. In particular, it is necessary to keep in mind that all the data we can collect contain imprecision. He wrote in \cite{Borel1909}

\begin{quote}
It is obvious that such a theory cannot correspond to any real and concrete problem; the conception of a denumerable sequence or unlimited sequence similar to the indefinite sequence of integers is a purely mathematical and theoretical conception, and the related speculations are themselves part of mathematics. It is only in an indirect way that these speculations, which constantly involve the notion of infinity, may be able to get practical applications\footnote{Il n'est pas besoin d'observer qu'une telle th\'eorie ne peut correspondre \`a aucun probl\`eme r\'eel et concret; la conception d'une suite d\'enombrable ou suite illimit\'ee analogue \`a a suite ind\'efinie des nombres entiers est une conception purement math\'ematique et th\'eorique, et les sp\'eculations qui s'y rattachent sont elles-m\^emes du domaine propre des math\'ematiques. C'est seulement d'une mani\`ere indirecte que ces sp\'eculations, o\`u intervient \`a chaque instant la notion de l'infini, pourront se r\'ev\'eler susceptibles d'applications pratiques.}.
\end{quote}

Among numerous other examples, let us also quote the following sentence from  his famous book `Le Hasard'  \cite{Borel1914} in which he commented on the solution of a functional equation describing a probability

\begin{quote}One should always accept with caution many results obtained in that way, even if the setting up of the functional equation cannot be criticized. One indeed is generally unaware of the extent to which a small error in a functional equation can imply considerable modifications in the results obtained with its solution\footnote{On doit toujours accepter avec circonspection beaucoup de r\'esultats obtenus de cette mani\`ere m\^eme si l'\'etablissement de l'\'equation fonctionnelle ne peut pas \^etre critiqu\'ee. On ignore en effet en g\'en\'eral \`a quel point une petite erreur peut entra\^\i ner des modifications consid\'erables dans les r\'esultats obtenus avec sa solution.}.\end{quote}

The strange kind of schizophrenia Borel developed towards probability had been the source of much comment. Borel was certainly the origin of a major renewal of the field of probability. He remained all his life a permanent supporter of the need for a probabilistic culture in society. But precisely for these reasons, the only possible justification for probability was for him its practical use. Therefore, the efforts needed to obtain a probability must be in direct connection with their practical importance. Because this importance must always be wisely considered,  Borel became reluctant to use high mathematics in any application, and first of all probability. This fact is maybe best illustrated by his own lack of consideration for the discipline - at least after WW1 when he told his wife that because high mathematics cracked his skull, he lazed about (in French {\it pantouflait}) in probability. Whether Camille Marbo's account was true or not, it is clear that Borel's viewpoint placed him in an awkward position with the new probabilistic generation. In particular, Borel's attitude profoundly shocked Paul L\'evy who became, as is well known, the major French specialist in probability of the inter-war period. For L\'evy,  mathematical constructions, based on sophisticated techniques if required,  were necessary to guarantee the validity of the application of probability theory. On this subject, consult \cite{BarbutLockerMazliak2004}. 

Considerations on randomness, probability and statistics were an important topic for papers from the beginning, especially Borel's in the {\it Revue du Mois}, the journal he founded with his wife in 1905. Borel found in it the exact tool he needed to improve his conceptions on these subjects. 
The first paper of the journal, written by Volterra\footnote{{\red About this article, we can consult the detailed study by Durand and Mazliak \cite{DurandMazliak2011}, which argues that the Prolusione was the Borel's first encounter with the kind of applied mathematics represented by biometrics. Volterra's quick affirmations about the role of probability and statistics in a near future had an influence on Borel's program. The fact that this article was the first published in the {\it Revue du Mois} gives it a strong programmatic nature.}} \cite{Volterra1906}, gave a hint of what Borel hoped to do: {\it For most mathematicians} wrote Volterra, {\it the desire emerges to lead their minds beyond the limits of pure mathematical analysis\footnote{Chez la plupart des math\'ematiciens s'\'eveille le d\'esir naturel de diriger leur esprit hors du cercle de la pure analyse math\'ematique.}}.  Let us observe that Volterra  quotes Pearson and the achievements of the British statistical school as an answer to a situation where the increasing complexity of biological knowledge makes Claude Bernard's experimental method on particular cases impractical. The biometric school provided statistical methods to overcome the necessary renunciation of {\it the hope of obtaining a mathematical relation between each ancestor and each descendent\footnote{L'esp\'erance d'obtenir une relation math\'ematique entre chaque anc\^etre et chaque descendant.}}. A hint of Borel's interest in Pearson's work is his publication in 1908 of a rather unnoticed note to the Comptes-Rendus in the section {\it Mathematical Statistics}, devoted to considerations bout Pearson's polymorphic curves (\cite{Borel1908b}).

\subsection{A social science : Discussions and Polemics}\label{polemics}

Borel's first intervention on probability in the {\it Revue du Mois} is \cite{Borel1906a}. The question of practical value was to remain at the core of Borel's philosophy of mathematical randomness until the end of his life. In practical life, Borel writes, probabilistic computations {\it put some elements of our decisions under a form easier to grasp. They replace some more or less complex data by a small number of simple figures\footnote{mettent sous une forme plus facile \`a saisir certains \'el\'ements de nos d\'ecisions. Ils remplacent certaines donn\'ees plus ou moins complexes par un petit nombre de chiffres simples.}.} 
Borel's main problem is therefore the scientific interpretation of statistical facts which are no longer considered only as a collection of descriptive data. 

In subsequent papers in the {\it Revue du Mois} (which became a major basis for the book {\it Le Hasard} \cite{Borel1914}), Borel was to develop these considerations. Let us in particular comment on two papers.

Paradoxes often appear in economic life. For instance, let us suppose that a milkman must decide on the retail price for one liter of milk if the wholesale price has been increased by 0.5 cent. If he decides to pass the cost along to his customers, he must necessarily amplify the increase of the retail price to at least one cent and may face a loss of customers. On the contrary he may accept a slight reduction of his profit until a subsequent increase of the wholesale price which would allow him to pass both increases to the retail price. There is no logical way to decide which choice is the best one.
In \cite{Borel1907}, Borel explained how a probability distribution obtained from statistical observation is the only satisfying answer. Borel connects the question with  Zeno's classical paradox of the heap of wheat : when do we decide that a certain number of grains of wheat constitute a heap? For Borel, again only a statistical answer is reasonable, through the interrogation of a sample of people. He wrote
\begin{quote}
The idea I would like to extract from what precedes, is that the adequate mathematical answer to numerous practical questions is a probability coefficient. Such an answer may not seem satisfactory to many minds who expect certainty from mathematics. It is a very unfortunate tendency : it is mostly regrettable that public education be, in that respect, so little developed; it is probably due to the fact that probability calculus is almost universally unknown, though it gradually penetrates everybody's life  (though various forms of insurance, mutual funds, pensions and so on). A probability coefficient constitutes a definitely  clear answer, corresponding to an absolutely tangible reality. Some minds will maintain that they `prefer' certainty; perhaps they also would prefer that 2 plus 2 were 5\footnote{L'id\'ee que je voudrais d\'egager de ce qui pr\'ec\`ede, c'est que la r\'eponse math\'ematique \`a donner \`a bien des questions pratiques est un coefficient de probabilit\'e. Une telle r\'eponse ne para\^\i tra pas satisfaisante \`a bien des esprits, qui attendent des math\'ematiques la certitude. C'est l\`a une tendance tr\`es f\^acheuse ; il est extr\^emement regrettable que l'\'education du public soit, \`a ce point de vue, si peu avanc\'ee; cela tient sans doute \`a ce que le calcul des probabilit\'es est \`a peu pr\`es universellement ignor\'e, bien qu'il p\'en\`etre chaque jour davantage dans la vie de chacun (assurances diverses, mutualit\'es, retraites, etc.). Un coefficient de probabilit\'e constitue une r\'eponse tout \`a fait claire, correspondant \`a une r\'ealit\'e absolument tangible. Certains esprits maintiendront qu'ils ``pr\'ef\`erent''  la certitude; ils pr\'ef\`ereraient peut-\^etre aussi que 2 et 2 fissent 5.}.
\end{quote}

In a subsequent paper (\cite{Borel1908}) Borel tried to demonstrate that better education on probability would be a good way to improve social solidarity and to avoid a {\it kind of individualism which is nothing but stupid selfishness}. For Borel, probability is the basis of what may be called {\it social mathematics}. One may state statistical laws satisfied by the facts, and such laws have the power to limit the excesses of individualism and therefore to contribute to the development of solidarity. 
For Borel quantification of randomness was a perfectly scientific attitude and he was probably surprised at stubborn reactions to this conception. He wrote:

\begin{quote}
Trying to go in depth into the reasons why the calculus of probability is unpleasant to many minds, I hope to succeed in showing that this unpleasantness is largely based on misunderstanding; I wish this misunderstanding  were removed, for the popularization of the conclusions, if not of the methods of this branch of science, could have much social utility. (\dots ) In general, men do not like to lose their names and to be represented by a number; neither do they like to be considered only a unit in a group without being individually identified. Here is already the reason why statistics are not popular and why easy jokes about them are generally welcomed. (\dots ) [The calculus of probability], not only lists past events, but also claims to predict future events to some extent; here lies its scientific aspect. This claim shocks the psychological sense of human freedom in the first place. (\dots ) One must not fear computations, if one is ready not to behave according to its results without having firstly weighted them for what they are : it is a strange illusion to think that individual independence is increased by ignorance. (\dots ) Someone who is starving is not much interested in the increase of average fortune : one must not seek in statistics or in computations arguments for the consolation of those who suffer from social inequalities; but this observation does not diminish anyhow the proper value of statistics or of the computations one uses to interpret them\footnote{En essayant d'approfondir les raisons pour lesquelles le calcul des probabilit\'es est antipathique \`a beaucoup d'esprits, j'esp\`ere arriver \`a faire voir que cette antipathie repose en grande partie sur un malentendu; il serait d\'esirable que ce malentendu soit dissip\'e, car la vulgarisation des conclusions, sinon des m\'ethodes de cette branche de la science, serait d'une grande utilit\'e sociale.
(\dots ) [L]'homme n'aime pas, en g\'en\'eral, perdre son nom et \^etre d\'esign\'e par un num\'ero; ni m\^eme \^etre compt\'e seulement comme une unit\'e dans un groupe sans \^etre individuellement d\'esign\'e. C'est d\'ej\`a l\`a une raison pour que la statistique ne soit pas populaire et pour que les plaisanteries faciles que l'on peut faire \`a son sujet soient g\'en\'eralement bien accueillies.
(\dots ) Le calcul des probabilit\'es  (\dots )  non content de recenser les \'ev\'enements pass\'es (\dots ) pr\'etend pr\'evoir dans une certaine mesure les \'ev\'enements futurs: c'est en cela qu'il est une science. Cette pr\'etention heurte tout d'abord le sentiment psychologique de la libert\'e humaine. (\dots )  On n'a rien \`a redouter du calcul, lorsqu'on est d\'ecid\'e \`a ne pas r\'egler sa conduite sur ses indications sans les avoir au pr\'ealable pes\'ees \`a leur juste valeur : c'est une illusion singuli\`ere de penser que l'ind\'ependance individuelle est accrue par l'ignorance. (\dots ) Celui qui meurt de faim s'int\'eresse peu \`a l'augmentation de la fortune moyenne: on ne doit pas chercher dans la statistique ni dans le calcul des arguments pour consoler ceux qui souffrent des in\'egalit\'es sociales; mais cette constatation ne diminue en rien la valeur propre des statistiques ni des calculs par lesquels on les interpr\`ete.}.
\end{quote}

We see how Borel followed Quetelet, Cournot and Markov's path in their quarrels with the opponents to the conception of an `average man', allegedly in contradiction to man's free will, whether a religious tenet or not : see \cite{Callens1997}, \cite{Martin2007} and \cite{Seneta2003} on this subject. 

{\red Borel was never really fond of abstract speculations about the meaning of probabilistic statements. In 1909 he published his {\it El\'ements de la théorie des probabilit\'es}, where he gave an exposition of the mathematical results forming the basis of his theory of probability. The book attracted Keynes's misunderstanding.  In two violent reviews (\cite{Keynes1910a} and \cite{Keynes1910b}), Keynes accused Borel of neglecting the philosophical interpretation of his results and of being concerned only by the mathematical side. Reading what Keynes wrote, it seems that Borel was for him  typical of the French way of considering these problems. Borel's answer came only fourteen years later (revenge is a cold dish!) when he strongly criticized Keynes  in  the review he wrote about the latter's 1921 treatise on Probability for the {\it Revue Philosophique} in 1924 \cite{Borel1924}. The amazingly harsh words used by Borel  in what is usually a rather conventional exercise are certainly a hint  that he was irritated by Keynes' obsession with the interpretation of results in probability at the expense of their use as scientific results (Keynes' aim in his treatise was to obtain a satisfactory basis for a logical interpretation of probability).  Borel was in particular upset that Keynes did not even mention statistical mechanics and the kinetic theory of gases which were for him, as we have seen, the paradigmatic examples which made randomness unavoidable in modern scientific conceptions. Borel wrote : {\it Maxwell is one of the most celebrated names the famous Cambridge university can be proud of, a university to which M.Keynes insists on connecting his intellectual education. If he is left aside, it is certainly not out of ignorance, but because systematically \footnote{Maxwell est cependant un des noms les plus illustres dont puisse s'enorgueillir la c\'el\`ebre Universit\'e de Cambridge , \`a laquelle M.Keynes tient \`a rattacher sa formation intellectuelle. S'il est laiss\'e de c\^ot\'e, ce n'est certainement pas par ignorance ni par oubli mais en vertu d'un syst\`eme.}} 
and later : {\it Is the highly successful application of probability to [physical] questions devoid of reality in M.Keynes' eyes?\footnote{Est-ce l'application pleine de succ\`es du calcul des probabilit\'es \`a ces questions qui est d\'epourvue de r\'ealit\'e aux yeux de M.Keynes?}}
Borel even went as far as considering Keynes' ideas as typical of the British way of thinking. As a conclusion to his introduction, Borel wrote the following lines:

\begin{quote}
This proves once again how different  British minds are from continental ones; we must neither be hypnotized by these differences nor look with obstinacy for what is to us incomprehensible; it is better to admit these differences as a {\it matter of fact\footnote{In English, in the original!}} and nevertheless try to adapt the original ideas of the British to our particular mentality. By doing so, we are almost certain to betray them, but at the same time we give to them the only chance to exert an influence on minds built differently from theirs. The history of science shows that this collaboration between minds which do not understand  each other perfectly is not only possible, but often fruitful. I strove to understand M.Keynes and, whenever I felt myself too far from him, I strove to faithfully look for the ideas, however different from his they may be, which had inspired me on reading his book\footnote{Ceci prouve une fois de plus combien sont diff\'erents les esprits des Anglais et les esprits des continentaux; nous ne devons pas nous hypnotiser sur ces diff\'erences et chercher avec obstination \`a comprendre ce qui est pour nous incompr\'ehensible; il vaut mieux admettre ces diff\'erences comme une {\it matter of fact} et essayer n\'eanmoins d'adapter \`a notre mentalit\'e particuli\`ere les id\'ees originales des Anglais. Ce faisant, nous sommes \`a peu pr\`es s\^urs de les trahir, mais en m\^eme temps de leur donner la seule chance qu'ils peuvent avoir d'exercer une infuence sur des esprits faits autrement que les leurs. L'histoire des sciences montre que cette collaboration entre esprits qui ne se comprennent pas compl\`etement est non seulement possible, mais souvent f\'econde. Je me suis efforc\'e de comprendre M.Keynes et, lorsque je me sentais trop \'eloign\'e de lui, de rechercher loyalement quelles id\'ees, peut-\^etre fort diff\'erentes des siennes, m'\'etaient sugg\'er\'ees par la lecture de son livre.
}.
\end{quote}}
As can be seen, Borel was not a man afraid of  controversies. On the contrary, the {\it Revue du Mois} had been the battlefield between Borel and several scholars of the time. One of the most important controversies was the recurrent discussion with biologist Le Dantec. As a brilliant mind, Le Dantec was a firm opponent of a scientific conception of probability.  Above all he objected to letting such a conception invade the domain of the life sciences. In 1910, Le Dantec wrote a paper called {\it Mathematicians and probability} \cite{LeDantec1910}\begin{quote}
I will teach everything necessary for actuaries, for the kinetic theory of gases, and so on, without having ever pronounced the dangerous word of probability or chance; I would rather call this part of mathematics : the computation of means, in the case of phenomena which are never submitted to any law\footnote{J'enseignerai ainsi tout ce qui est n\'ecessaire pour les actuaires, pour la th\'eorie cin\'etique des gaz, etc., sans avoir jamais prononc\'e le mot dangereux de probabilit\'e ou de chance; j'appellerais volontiers cette partie des math\'ematiques : le calcul des moyennes, dans le cas des ph\'enom\`enes qui ne sont jamais soumis \`a aucune loi}
\end{quote}

We shall not comment further on Le Dantec's discussion with Borel which is thoroughly analyzed in \cite{BruBruChung1999}. 

Another discussion was Borel's direct opposition to Bergson's conception of education. This theme was examined by Callens in \cite{Callens1997}. We will just recall some of its basic features. The French philosopher had adopted a characteristic anti-intellectual attitude and defended the idea that an {\it intellectual man} deprives himself of the natural intelligence which alone can provide tools for an adaptation to the environment. An improvement of this natural intelligence can only be achieved through classical studies (Greek, Latin, geometry\dots ) which allow the selection of an elite of men of action. As Grivet wrote when he commented on Bergson's philosophy of action : { \it Forget about adjusting your action to knowledge, or you will condemn yourself  to never do anything. The man of action is the one who, at a given time, knows how to silence the faculty of reasoning.\footnote{Renoncez \`a r\'egler votre action sur la connaissance, ou vous vous exposerez \`a ne jamais rien faire. L'homme d'action est celui qui, \`a un moment donn\'e, sait faire taire au plus vite la facult\'e de raisonner (\cite{Grivet1911}, pp.471-472).}}
Obviously, quantifying randomness is quite unnatural to the human mind and therefore does not fit Bergson's conceptions. For Borel, on the contrary, it is precisely from data mining and statistical treatment that the good sailor can prepare for the storm. The practical value of probability for Borel is that computation allows one not to bury one's head in the sand. It is therefore the best reliable support for a clever and brave attitude in face of danger. Bergson advocates that one should follow one's instinct. But, Borel asserts, instinct deceives - and such a behavior would lead to an attitude much less adapted to  life than being trained to risk. Borel finds a proof of this deceiving aspect of instinct in the way men eventually had to renounce a completely mechanical explanation of the universe. For example, {\it [t]he discovery and study of radioactivity showed that mechanical explanations are sometimes certainly insufficient and must 
give way to statistical explanations\footnote{L'explication m\'ecanique de l'univers s'est toujours d\'erob\'ee jusqu'au jour o\`u la d\'ecouverte et l'\'etude de la radioactivit\'e ont montr\'e que les explications m\'ecaniques sont parfois certainement insuffisantes et doivent alors c\'eder le pas aux explications statistiques (\cite{Borel1914},p.iii)}}.
Much later, when Borel wrote the final part of the long series of volumes on probability theory he edited, which he devoted to the {\it Practical and philosophical value of probability}, he summed up the views he had developed in his previous publications. {\it Probability calculus is essentially social science\footnote{Le calcul des probabilit\'es est essentiellement une science sociale (\cite{Borel1939}, p.129)}}
which allows the simulation of the social effects of a decision. This is profoundly related to the fact that probability uses sets and populations and not individuals as basic components. Thus {\it the theory of probability dominates all experimental science as deductive logic dominates mathematical science \footnote{la th\'eorie des probabilit\'es domine toute la science exp\'erimentale, autant que la logique d\'eductive domine la science math\'ematique (\cite{Borel1939}, p.126)}}. 
Probability is in fact the only way to give a meaning to a collection of data. Borel is extremely assertive in his conclusion : {\it The value of any science is founded on reasoning with probabilities.\footnote{La valeur de toute science a comme fondement des raisonnements de probabilit\'es \cite{Borel1939}, p.126}}. 

{\red As we have seen, in 1905 Borel began his interest in probability in a purely mathematical way - the application of measure theory. But, as the practical and social value of mathematics had always been a concern for him,  questions about justifying the application of probability theory to statistics - which deals with the concrete aspect of random phenomena - naturally came soon. An  exchange with Lucien March (the head of the Statistique G\'en\'erale de France) in 1907\footnote{Letter from March to Borel, 16 July 1907. Archive of Paris Academy of Sciences.} testified to this fact. March explained that the core of the application is the assimilation between observed frequencies and probabilities. This could certainly be sometimes true but the problem had not been studied enough. In the next years, Borel was led to see some probabilistic results as the theoretical basis for the study of statistical situations. Probability became therefore in Borel's mind this ``father of mathematical statistics`` (p\'ere de la statistique math\'ematique), described later by Fr\'echet in 1948 (\cite{Frechet1948}). Borel was to advocate this idea in his own 1923 presidential speech. He provided  the law of errors as main example,  the use of which he considered universal in statistics (\cite{Borel1923}). Not coincidentally, in the same year,  Borel and Deltheil  wrote the little book {\it Probabilit\'es, Erreurs} (\cite{BorelDeltheil1924}) where they aimed at presenting the Central Limit Theorem and some of its applications with the least  mathematical technique.}

\subsection{The war experience}

The sudden outbreak of the Great War was the major event which allowed Borel to make a kind of full-sized test of his former considerations on randomness and statistics. In a paper \cite{Borel1919} written in 1919 for the revival of the {\it Revue du Mois} \footnote{The publication of the journal in fact lasted only one year, mostly because the energy and the faith of the beginning had vanished.}, Borel wrote that {\it the Great War has been a decisive experiment as well as an education}\footnote{la grande guerre a \'et\'e une \'epreuve d\'ecisive et en m\^eme temps un enseignement} to help us realize that science had a major role in the material development of mankind, an idea commonly accepted during the 19th Century but mostly limited to intellectual considerations that {\it had not deeply penetrated the intimate depths of  consciousness}\footnote{n'avait pas p\'en\'etr\'e dans les profondeurs intimes de la conscience}. 

The recently edited war-correspondence between Volterra and Borel (\cite{MazliakTazzioli2009}) provides an impressive picture of how, since the very beginning, Borel had decided to involve himself in the war effort, though at the age of 43 he could have made a more comfortable choice. But in fact, as some comments in \cite{MazliakTazzioli2009} try to illustrate, Borel's conception of social life prevented him from even thinking of the possibility of staying on the sidelines. Just as the mathematician needed to be implied in the life of the city (and Borel's leitmotiv was that probability theory gave a tool for that to mathematicians), the scientist could not help playing his part in the war. 
Borel enlisted in the army in 1915 and participated in the testing of sound ranging apparatus, in particular the so-called Cotton-Weiss apparatus devised in 1914 by physicists Aim\'e Cotton and Pierre Weiss for the localization of artillery batteries, on the front.

But probably more significant to the story we present here, Borel was soon proposed to take up important responsibilities at a high governmental level, close to the center of power in the very centralized France of that time. In November 1915, his friend and colleague mathematician, Paul Painlev\'e, minister for public instruction, decided to set up a special service connected to his ministry in order to gather the various enterprises dealing with technical research into the war effort. Painlev\'e and Borel were very close friends, belonging to the same political family  and sharing common views on the role of the engaged scientist in the city. Painlev\'e contributed several times to the {\it Revue du Mois}; in 1913, together with Borel and Charles Maurain, he wrote a visionary book on the development of aircraft \cite{BorelMaurainPainleve1913}. Painlev\'e asked Borel to head the special service called the {\it Direction des inventions int\'eressant la d\'efense nationale} (Direction of Inventions related to national defense).  A clipping from the {\it Journal de Paris} on 15 November 1915 asserts that  Borel's nomination
\begin{quote} as head of the service of inventions will be welcomed with general satisfaction. He is known in the scientific world for his works of pure mathematics and mechanics. Also, he has already collaborated with the present minister for public education in studies about aviation\footnote{Sa nomination \`a la t\^ete du service des inventions sera accueillie avec une satisfaction unanime. Il est connu dans le monde scientifique pour ses travaux de math\'ematique pure et de m\'ecanique. Enfin, il a collabor\'e d\'ej\`a avec le ministre actuel de l'instruction publique, dans les \'etudes que celui-ci a faites sur l'aviation}.
\end{quote}

The journalist obviously could not have been aware of the importance of Borel's reflections on statistics and randomness though they were probably constantly on Borel's mind. {\red The same year Borel was chosen to be a member of the SSP board. His new function as head of the {\it Direction des Inventions} was probably the main reason for this choice, but Borel's interest in statistics was certainly taken into consideration.} When, two years later, Painlev\'e became Prime Minister ({\it Pr\'esident du Conseil} - President of the council, during the 3rd Republic) and Borel was promoted as head of Cabinet, he assigned himself a major task : to organize rationally the statistical data mining devoted to helping the Prime Minister in his decisions. In fact, Borel did not have enough time for that. Three months after his arrival as Prime Minister, in November 1917, Painlev\'e's government was toppled  as a consequence of the Italian disaster of Caporetto (see \cite{MazliakTazzioli2009} for details), and Borel returned to his experiments on sound ranging on the front. However, this experience made a profound impression on him and just after the war, he wrote a text, published in the Journal of Paris Statistical Society in January 1920, entitled {\it Statistics and the Organization of the Presidency of the Council}(\cite{Borel1920}). 
In this text, Borel wrote the following lines directly inspired from his observations when working with Painlev\'e. 

\begin{quote}
I would like to insist on the role of the service depending on the Presidency of Council which we can call, to specify its nature, a statistical cabinet, as it is for me  this statistical cabinet that befalls one of the most important and at the same time one of the most delicate tasks of the government of the country. The number and the material importance of statistical documents increase each day in every country; one  realizes better, indeed, the importance of obtaining sufficiently detailed statistics in order to use them in different situations. Social phenomena are too complex to be easily included in oversimplified formulae. But, on the other hand, reading and interpreting considerable statistical documents require not only some specific form of education but also a large amount of time. We must admit that the heads of Government are not short of specialized education but they are mainly short of time. It is therefore necessary that men in which they have full trust sum up and interpret statistical documents  for them.  But, once data are summed up and interpreted, any rigorously scientific and objective work is impossible; it is therefore not possible to entrust civil servants with this work whatever their professional value can be, as their opinions may be on such or such question of economic, customs or fiscal  politics, opposed to the Government's opinion. (\dots ) Here is not the place to insist on the fact that statistics is an indispensable aid to all those who are in charge of the heavy task of ruling a country. If nevertheless we all agree on the principle, there may be divergences of opinion on the most favorable modes of carrying things out,  that may lead to a very profitable discussion. The times seem to me particularly well suited for this discussion   as within some weeks, France will choose the Government  which will organize peace after choosing the Governments which have won the war. Without involving ourselves here in political matters, we can however assert that whatever the Government will be, a clever use of statistics will be useful to them. (\dots ) Let me remind you, as an example of what can be done in that direction, of the organization of the technical services of the  War cabinet such as Painlev\'e  conceived it. The numerous statistical documents related to the war policies (French, allies or enemy troop sizes, losses, ammunition, necessary and available tonnages, submarine warfare, exchange rates and so on) were collected and summed up. A notebook, called the black notebook, was constituted; it was formed of around ten cardboard sheets, equipped with tabs and also a folder. By opening the notebook on the page of troop sizes for instance, one immediately found the current, summarized information in the specific form required by the minister on a sheet of paper updated  each week; in the corresponding folder, were collected retrospective information, further details, graphs.  (\dots ) When M.Painlev\'e became also Head of the Government, a grey notebook, of the same kind and containing statistical documents of interest for the economic committee, was joined to the black one which contained the statistics of interest for the war committee\footnote{Je voudrais (\dots ) insister un peu sur le r\^ole de l'organe de la pr\'esidence du Conseil que nous pouvons appeler, pour pr\'eciser sa nature, cabinet statistique, car c'est, \`a mon avis, \`a ce cabinet statistique que doit incomber une des t\^aches les plus importantes et en m\^eme temps les plus d\'elicates dans le gouvernement du pays. Le nombre et l'importance mat\'erielle des documents statistiques augmente chaque jour dans tous les pays; on se rend mieux compte, en effet, de l'importance qu'il y a \`a poss\'eder des statistiques suffisamment d\'etaill\'ees pour qu'elles soient utilisables \`a des fins diverses. Les ph\'enom\`enes sociaux sont trop complexes pour qu'il soit possible de les enfermer dans des formules trop simplifi\'ees. Mais d'autre part, pour lire et interpr\'eter des documents statistiques consid\'erables, il faut, non seulement une \'education sp\'eciale, mais beaucoup de temps. Nous devons admettre que l'\'education sp\'eciale ne fait pas d\'efaut aux chefs du Gouvernement, mais c'est le temps qui leur manque le plus. Il est donc n\'ecessaire que les hommes en qui ils aient pleine confiance r\'esument et interpr\`etent pour eux les documents statistiques. Or, d\`es qu'il y a r\'esum\'e et interpr\'etation, il ne peut plus \^etre question d'un travail rigoureusement scientifique et objectif; il n'est donc pas possible de confier ce travail \`a des fonctionnaires quelle que soit leur valeur professionnelle, dont les vues personnelles peuvent \^etre, sur telle question de politique \'economique, douani\`ere ou fiscale, en opposition avec celles de Gouvernement. (\dots )
[C]e n'est pas ici qu'il est n\'ecessaire d'insister sur le fait que la statistique est un auxilliaire indispensable pour ceux qui assument la lourde t\^ache de gouverner un pays. Si cependant nous sommes tous d'accord sur le principe, il peut y avoir sur les modes d'ex\'ecution les plus favorables, des divergences d'appr\'eciation qui pourraient conduire \`a une discussion tr\`es profitable. Le moment me para\^\i t particuil\`erement bien choisi pour cette discussion car c'est dans quelques semaines que la France va, apr\`es les Gouvernements qui ont gagn\'e la guerre, conna\^\i tre les Gouvernements qui organiseront la paix. Nous n'avons pas \`a intervenir ici dans les questions politiques, mais nous pouvons affirmer que, quels que soient ces Gouvernements, l'emploi judicieux des statistiques leur sera n\'ecessaire. (\dots ) Qu'il me soit permis de rappeler, comme exemple de ce qui peut \^etre fait dans ce sens, l'organisation des services techniques du Cabinet au minist\`ere de la Guerre, telle que l'avait con\c c ue M.Painlev\'e. Les tr\`es nombreux documents statistiques int\'eressant la politique de guerre (effectifs fran\c c ais, alli\'es ou ennemis, pertes, munitions, tonnages n\'ecessaires et disponibles, guerre sous-marine, changes, etc. ) \'etaient rassembl\'es et r\'esum\'es. Un cahier, que l'on appelait le cahier noir, avait \'et\'e constitu\'e; il \'etait form\'e d'une dizaine de feuilles de carton, munies d'onglets et dont chacune comportait en outre une pochette. En ouvrant le cahier \`a la page des effectifs, par exemple, on trouvait imm\'ediatement sur une feuille renouvel\'ee chaque semaine, les renseignements actuels r\'esum\'es sous la forme d\'esir\'ee par le ministre; dans la pochette correspondante se trouvaient des renseignements r\'etrospectifs, des d\'etails compl\'ementaires, des graphiques. (\dots ) Lorsque M.Painlev\'e joignit la pr\'esidence du Conseil au minist\`ere de la Guerre, au cahier noir qui contenait les statistiques int\'eressant le Comit\'e de Guerre fut adjoint un "cahier gris" \'etabli sur le m\^eme mod\`ele et renfermant les documents statistiques int\'eressant le Comit\'e \'economique.}.
\end{quote}

After he became a politician in the classical meaning of the word (first as Mayor of St Affrique in 1923, then as Deputy in 1924 - he was even a minister for the Navy in Painlev\'e's (once again short) presidency of the Council in 1925), Borel tried to obtain the creation of a service of statistical documentation and economic studies. He pleaded therefore for the {\it institution of data, given that the 19th Century, in its romanticism, had only established an institution of facts \footnote{ pour l'institution des donn\'ees, alors que le XIX\`eme si\`ecle, par romantisme, n'avait r\'ealis\'e qu'une institution des faits.}}. In his posthumous comments on Borel's works, Denjoy  wrote that he doubted {\it that a single analyst ever had a keener sense of numerical reality than Borel}
\footnote{Je doute qu'un seul analyste ait eu autant que Borel, ni m\^eme \`a un degr\'e comparable au sien le sens de la r\'ealit\'e num\'erique}(\cite{Denjoy1972}). 

The war experience had been for Borel a large-scale numerical experience. This was not only because total warfare created a situation where enormous quantities of materiel, food, and weapons were used, replaced, and exchanged,  and in which the unprecedented size of the armies made it necessary to develop tools for the direction of millions of soldiers; but also because the trauma caused by the huge losses forced the forging of a new conception of social life. We refer the reader to \cite{MazliakTazzioli2009} and \cite{Mazliak2009b} 
for a presentation (and some affecting documents) on how Borel faced the tragedy - a tragedy which besides affected him personally through the loss of his foster son Fernand as well as it affected his fellow mathematicians Hadamard, Picard and so many others. For Borel, social mathematics belonged to the tools provided to society to avoid the woes created by man's natural selfishness. It remains profoundly linked to a constitutional vision of politics, close to the ideals promoted by the social-radical politicians of his generation, providing a tool whose purpose was among other things to keep any authoritarian tendency at bay. In a book written in 1925, whose title {\it Organize} is in itself a program, Borel advocated an international law to which each State must be submitted. Each State must be accountable to public opinion and to an assembly of nations. A combination of exchanges between {\it assemblies of free men}, seen as a combination of economic, cultural, and political exchanges is seen as the best guarantee against the {\it most absolute powers}. Some letters exchanged with Volterra  (who entered very soon into direct opposition to  Mussolini's regime - see \cite{Goodman2007} and \cite{GuerraggioPaoloni2008}) bear testimony of how the Italian events  at the beginning of the 1920s were for Borel an illustration of the aforementioned necessity. We shall see in the next section how the foundation of the Henri Poincar\'e Institute was used by Borel as an opportunity to invite fellow scientists trapped in the convulsions of the interwar period in Europe.

\section{The emergence of mathematical statistics in France}

In 1920, Borel, who occupied the chair of function theory asked the University of Paris to be transferred to the chair of Probability Calculus and Mathematical Physics. This happened at the precise moment when he decided to resign from the position of vice-director of the Ecole Normale Sup\'erieure, because, wrote Marbo in her memoirs (\cite{Marbo1967}), he could not any longer face  the ghosts of all the young students fallen during the war (see \cite{Mazliak2009}). We have already mentioned how Marbo commented on Borel's weariness about higher mathematics after the war. Several commentators take this assertion for granted and want to see Borel's desire for transfer as a proof of this lack of interest, due to his ambivalent relationship towards probability as a mathematical theory. Nevertheless, in a letter to Volterra (\cite{MazliakTazzioli2009}, p.138 ), Borel told his Italian colleague that he had tried to convince Langevin to accept the position. Had the latter accepted, it would be interesting to understand what Borel would have done concerning his own situation. But Langevin had refused and Borel asked for the position which he also probably felt as a convenient place to develop his ideas on randomness. Gispert (\cite{Gispert2009}) and Gispert and Leloup \cite{GispertLeloup2009} have observed that the discussion about Borel's transfer was the occasion of a new crisis between Borel and Lebesgue after their break in 1917, Borel having proposed the transformation of his chair of theory of functions into a chair of theoretical and celestial physics. Lebesgue was upset because he saw this as a betrayal of mathematics. He wrote : {\it Should M.Borel succeed - and he will - his lectures will attract physicists and if he attracts also mathematicians, they will be lost for mathematics}, illustrating the poor opinion he himself had of probability.

Borel was intent in doing all that was possible to improve the diffusion of the mathematics of randomness in France. A trace of these efforts can be found in the introduction to the small textbook he wrote with Deltheil (\cite{BorelDeltheil1924}). 

\begin{quote}
It is only for reasons of tradition, shall we say of routine, that the elements of probability calculus are not present in the syllabus of secondary schools, where they could with benefit replace many of the  remaining subjects, which are still there only because no one thought of getting rid of them\footnote{C'est uniquement pour des raisons de tradition - l'on n'ose \'ecrire de routine - que les \'el\'ements du calcul des probabilit\'es ne figurent pas au programme de l'enseignement secondaire, o\`u ils remplaceraient avantageusement bien des mati\`eres qui y subsistent pour le seul motif que personne ne se donne la peine de les supprimer.}.
\end{quote}\rm

\subsection{Strasbourg}

An important lesson of the war had been the observation of the efficiency of German organization. During the entire 19th Century, and especially after the  birth of the German Empire, the development of Germany had been accompanied by a development of powerful statistical institutions (see \cite{Desrosieres2000}, p.218 to 231). As we have seen, Borel also knew about the active school of biometrics in Great Britain around Karl Pearson. By contrast, in France, the scene of academic mathematical statistics was desperately empty. 

Just after the war, the discipline experienced a promising development in Strasbourg, which had just been returned to France. 

In 1919, Fr\'echet was sent to the University of Strasbourg, and came there as a missionary of science. The Government wanted the institution to be a showcase of French research success. Deputy Manoury wrote the following letter on 5 April 1919 to the French government representative Alexandre Millerand:

\begin{quote}
You know better than anyone the considerable importance given by the Germans to this university and the special attention they paid to making it one of the most brilliant, if the not the brightest, of all the universities in the Empire. You have certainly also seen that they predicted that in less than 3 years France would have sabotaged their work. How could we face this challenge ?\footnote{Vous savez mieux que personne l'importance consid\'erable que les allemands avaient donn\'ee \`a cette universit\'e et la coquetterie qu'ils ont mise \`a en faire une des plus brillantes sinon la plus brillante de l'empire. Vous avez certainement vu aussi qu'ils ont pr\'edit en partant qu'en moins de 3 ans la France aurait sabot\'e leur \oe uvre. Comment relever ce d\'efi ?}
\end{quote}

On the university of Strasbourg at that time,  we confer the reader to the book \cite{Strasbourg2005}, and more precisely to Siegmund-Schultze's text about Fr\'echet \cite{SiegmundSchultze2005}. 

Thus, Strasbourg University had become within 10 years a first-rank institution in France, and a place of original intellectual experiments. During the Imperial period, Stra\ss burg had indeed been an important place for statistical research, with Lexis and Knapp. The young statistician Henri Bunle was sent to Strasbourg to retrieve German technical knowledge, as he explained in an exciting interview in 1982 (\cite{Desrosieres2005}). Desrosi\`eres asked Bunle what his task in Strasbourg was in 1919, and Bunle explained that his mission was to control the German statistical bureau which published a statistical directory for Alsace-Lorraine where 7 or 8 people worked. Bunle said: 

\begin{quote}
I recruited  Alsatians and Lothringians with deep [local] roots. I put these guys next to [the Germans]. I  went to see the Germans and said : I have placed people from Alsace and Lorraine beside you so that you can explain to them thoroughly what you are doing. You will leave only when these people tell me that they know the trade. Thus, as they wanted to go away, everything turned out OK.\footnote{ J'ai recrut\'e des Alsaciens -Lorrains de bonne souche. Je leur ai mis des types \`a c\^ot\'e. J'ai \'et\'e voir les Allemands et je leur ai dit : je vous ai mis des Alsacien-Lorrains \`a c\^ot\'e de vous pour que vous les mettiez enti\`erement au courant  ce que vous avez \`a faire. Vous ne partirez que lorsque ces gens me diront qu'ils connaissent le m\'etier. Alors, comme ils voulaient s'en aller, \c ca s'est bien pass\'e.}.
\end{quote}

Among the new pedagogical initiatives in Strasbourg, there was the creation of the {\it Institut d'\'etudes commerciales } (Institute for Commercial Studies) where Fr\'echet and the sociologist Maurice Halbwachs taught, as soon as 1920. Later, they published a book on their common experience (\cite{FrechetHalbwachs1924}). In the preface they explain the purpose of the book. For Fr\'echet, who was in charge of the lectures on insurance, the scientist who is deeply involved in speculative research must not lose interest in practice, and it is useful for the progress of science to spread its results. For Halbwachs, the statistical method is only a routine for the one who cannot catch its spirit and its deep scientific meaning. Halbwachs, though he had been a student in literature, had extensively reflected on the meaning of statistical process in a study devoted to social life, in particular in his book on Quetelet published at the eve of the war. The idea of Fr\'echet and  Halbwachs' book is to present the principles of probability and their application using only the most basic notions of mathematics.

{\red It is therefore interesting to observe that Fr\'echet and Halbwachs published a book on probability as a result of their teaching of statistics at Strasbourg institute. Fr\'echet accepted from the beginning that mathematical statistics was generated by probability, as he claimed later in his 1948 address (\cite{Frechet1948}). Even if he began to concentrate on statistics during his Strasbourg period, Fr\'echet was essentially concerned by probability theory before the 1930s\footnote{See \cite{HavlovaMazliakSisma2005} for details on Fr\'echet's beginnings on probability and in particular his  correspondence with Hostinsk\'y about Markov chains.}. His statistical activity is therefore mostly situated after the period of time we consider in the present paper, as is proven by his  late election at the head of the SSP at the venerable age of 70. On Fr\'echet's statistical works, consult \cite{Armatte2001} and \cite{Barbut2007}. }

\subsection{Teaching statistics. The ISUP}

 Borel was not insensitive to Strasbourg's innovations. He used his new political and scientific influence (he had been elected to the Academy of Sciences in 1921)  to help promote the teaching of mathematical statistics in Paris. As a result  the {\it Institut de Statistique de l'Universit\'e de Paris}  (ISUP - Statistical Institute of Paris University) was created in 1922. 

Emile Borel convinced Lucien March and Fernand Faure from the Faculty of Law to join him in setting up the new institution.  The Institute depended on the four faculties of Science, Medicine, Law and Literature. It was seen as a promising example of an interdisciplinary place for teaching. Such an institute perfectly illustrated the idea of knowledge, exchange and meeting - especially in the field of probability and statistics - that Borel sought to consolidate in France. The creation of the ISUP remained a well-kept secret, though it was mentioned for example in the journal {\it Vient de para\^\i tre} in January 1923, where a particular emphasis was put on the innovation it represented.

\begin{quote}
Paris University has just created a Statistical  Institute, where  the Faculty of Law and the Faculty of Science will collaborate; its head office will be in the buildings of the Faculty of Law. At the same time, the construction of the Institute of Geography on rue Pierre Curie,  where  the Faculty of Science and the Faculty of Literature will fraternize, is being completed. When numerous institutes of that kind are truly alive, the word University will no more be only a word, meaning the reunion of several Faculties ignoring each other. Several Universities outside Paris have been for a long time, on that point as well as on others, ahead of the University of Paris\footnote{L'Universit\'e de Paris vient de cr\'eer un Institut de Statistique, o\`u collaboreront la Facult\'e de Droit et la Facult\'e des Sciences ; son si\`ege sera dans les b\^atiments de la Facult\'e de Droit. En m\^eme temps, s'ach\`event rue Pierre Curie, les b\^atiments de l'Institut de G\'eographie o\`u fraterniseront la Facult\'e des Sciences et la Facult\'e des lettres. Lorsque de nombreux Instituts de ce genre seront bien vivants, le mot d'Universit\'e ne sera plus un simple mot, d\'esignant la r\'eunion de plusieurs Facult\'es s'ignorant les unes les autres. Certaines Universit\'es des d\'epartements ont, sur ce point comme sur d'autres, devanc\'e depuis longtemps l'Universit\'e de Paris.}. 
\end{quote}

In his paper \cite{Meusnier2006}  Meusnier observed that at the ISUP the students were trained in four main subjects : demography and economy, actuarial sciences, industrial technique and research, and medicine. The aim of the institute was naturally to teach statistics, with officially both a theoretical and practical point of view. Nevertheless, the program during the first years of its existence showed a clear desire to insist on applications. The latter concern a wide range of domains such as finance, political economy, demography, biometrics, public health, prediction, insurance, trading, agriculture, transports, bank, credit or public finances.

Although the ISUP officially opened in 1922, the full program really began in 1924-1925. In the two previous years, there had been some lectures, but no diplomas were delivered. At the beginning, Borel taught the course of \textit{statistical methods}, which dealt with mathematical statistics, 

The number of students at the ISUP during the first years remained extremely low. Only 4 students were registered in the academic year 1924-25,  when the first  diploma was delivered at the Institute. The fact probably fostered the skepticism with which the creation of the Institute was met by field statisticians - a skepticism obvious in the already mentioned interview of Henri Bunle \cite{Desrosieres2005}; though the causality could also be in the other direction.

 \begin{quote}
 With M.Borel, M.March has created the statistical whatchamacallit of the Sorbonne. M.Borel read a small course during one year, and then M.March, and then M.Huber on demography. When Borel was fed up, he passed on his course to Darmois who was in Nancy. Darmois began to collect information about what had been done in the United Kingdom. Because in the United Kingdom, they had been working more. There was a book for teaching statistics. He began to teach statistics. Besides, there is a book by him. He developed the mathematical point of view a bit more later. This was the ISUP regime\footnote{Avec Monsieur Borel, M.March a cr\'e\'e le machin de Statistique de la Sorbonne. Monsieur Borel a fait un petit cours pendant une ann\'ee, et puis M.March, puis M. Huber sur la d\'emographie. Quand Borel en a eu assez, il a pass\'e son cours \`a Darmois qui \'etait \`a Nancy. Darmois a commenc\'e \`a se mettre au courant de ce qui avait commenc\'e \`a se faire en Angleterre. Parce qu'en Angleterre, ils avaient travaill\'e davantage. Il y avait un volume pour l'enseignement de la Statistique. Il a commenc\'e \`a enseigner la Statistique. D'ailleurs, il y a un volume de lui. Il a d\'evelopp\'e un peu plus du point de vue math\'ematique. Voil\`a le r\'egime de l'ISUP.}.
\end{quote}

As Bunle mentioned, already in 1924 Borel asked Darmois to replace him.  

\subsection{Darmois}
George Darmois was born on 24 June 1888 in Eply near Nancy. In 1906,  he entered the \'Ecole Normale Sup\'erieure. Under Darboux's influence, in 1911 he began a doctorate on a subject mixing geometry and analysis. It concerned the study of partial differential equations arising from geometry. The title of his thesis was \textit{About algebraic curves with constant torsion}. However, due to the First World War, Darmois defended his thesis only on  26 February 1921. Nevertheless, the War had other, more important influences on his scientific evolution. During the War, Darmois had worked on two military problems, one about sound ranging and the other about ballistics. He described how these studies had played a decisive role (\cite{Darmois1955}):

\begin{quote}
The war of 1914-1918, having oriented me towards ballistics and shooting problems, and then towards location by sound and the problems of measuring and of wave propagation, had deeply inflected my spirit towards mathematical physics and the calculus of probability\footnote{La guerre de 1914-1918, en m'orientant vers la balistique et les probl\`emes de tir, puis vers le rep\'erage par le son et les probl\`emes de mesure et de propagation des ondes, a tr\`es fortement infl\'echi mon esprit vers la Physique math\'ematique et le Calcul des probabilit\'es}.
\end{quote} 

After the war, Darmois became a professor of analysis in Nancy's University of Sciences.  During the 1920s, Darmois' primary interest was mathematical physics, especially the theory of relativity. He published several articles about relativity, and as lately as 1930, he made a conference where he spoke about the experimental verifications of this theory and published a book in 1932 on the subject \cite{Darmois1932b}. The peculiar destiny of relativity theory in France during the 1920s was studied in \cite{Ritter2009} through Eyraud's case. Eyraud was another young mathematician  who subsequently turned towards statistics and founded a statistical institute in Lyon in 1930, the {\it Institut de Science Financi\`ere et d'Assurances} (Institute for Financial Science and Insurances), the first French institution to deliver a diploma for actuaries. In fact, Darmois was interested in any experimental aspect of science, not only in physics.  In 1923, he was already teaching probabilities and  applications to statistics. Darmois later explained (\cite{Darmois1955})

\begin{quote}
The decision I took in Nancy in 1923 to connect the teaching and investigations of probability calculus with several applications to statistics stemmed from the desire of constituting in France a school of theoretical and practical statistics. Great Britain and the United States showed the way, and it was important to follow their examples. I was thus led to participate in national, and later international, statistical activities\footnote{La d\'ecision que j'ai prise \`a Nancy vers 1923 de joindre aux enseignements et recherches sur le calcul des probabilit\'es diverses applications \`a la statistique, est venue du d\'esir de constituer en France une \'ecole de Statistique th\'eorique et pratique. La Grande-Bretagne et les \'Etats-Unis montraient le chemin, il importait de suivre leur exemple. J'ai ainsi \'et\'e amen\'e \`a prendre part aux activit\'es statistiques nationales, puis internationales.}.
\end{quote}

Apart from Darmois' intellectual interests for statistics, another reason why Borel called Darmois to Paris may have been the ties he had with the industrial world. In a recently published interview, Guilbaud informs us (\cite{Guilbaud2008}) that Darmois owned a small forge in the Vosges and frequently mocked his colleagues who disdained contact with industrialists.
 
 In 1928, at the International Congress of Mathematicians in Bologna, Darmois gave a talk entitled \textit{About the analysis and comparison of statistical series which could be developed in time (\textit{the time correlation problem})}.
 As Danjon wrote (\cite{Danjon1960})

\begin{quote}
Georges Darmois set himself two tasks. Firstly, to disseminate the power of statistical methods applied to sciences of observation, to biometrics, to applied psychology, to econometrics, to production control,  to operational research, and so on. The astronomers cannot forget that he was the first in France to teach statistics and stellar dynamics as early as 1928-29. His apostolic mission for statistics was pursued relentlessly in the form of lectures, conferences, seminars. In some domains, success was instantaneous; George Darmois' ideas,  by dint of obstinacy, managed to permeate the least prepared environments and eventually made them prevail against routine. Moreover, he tried to improve those parts of the theory which seemed to deserve special attention. And so he especially cared about the general theories of estimation from random sampling, of which the general theory of errors is a particular case\footnote{George Darmois s'\'etait assign\'e deux t\^aches. En premier lieu, faire comma\^\i tre la puissance des m\'ethodes statistiques appliqu\'ees aux Sciences d'observation, \`a la Biom\'etrie, \`a la Psychologie Aplliqu\'ee, \`a l'\'Econom\'etrie, au contr\^ole des entreprises, \`a la Recherche op\'erationnelle, etc. Les astronomes ne peuvent oublier qu'il fut le premier en France \`a professer la Statistiques et la Dynamique stellaire, d\`es 1928-1929. Son apostolat en faveur de la Statistique s'est pourquivi sans rel\^ache, sous la forme de cours, de conf\'erences, de s\'eminaires. Dans certains domaines, le succ\`es en fut imm\'ediat ; l'obstination de Geroge Darmois devait faire p\'en\'etrer ses id\'eees dans les milieux les moins bien pr\'epar\'es et finalement les faire pr\'evaloir contre la routine.Il s'est \'efforc\'e, en outre, de faire progresser les parties de la th\'eorie qui lui paraissaient m\'eriter une attention sp\'eciale. C'est ainsi qu'il s'est particuli\`erement appliqu\'e aux th\'eories g\'en\'erales de l'estimation sur \'echantillon al\'eatoire, dont la th\'eorie g\'en\'erale des erreurs est un cas particulier.}.
\end{quote}

  In 1925-1926, the ISUP syllabus was clearly oriented towards applied statistics and the probabilistic content was rather elementary and superficial. Meusnier (\cite{Meusnier2006}) listed the different items of the syllabus, which appear more as a collection of independent tools than as a real continuous exposition of a theory. It is worth noticing that the definition and {\it fundamental principles} of probability appear only as the 7th item of the syllabus. Moreover, the only theorem of probability theory mentioned is the (weak) law of large numbers. 

At the ISUP,  Darmois' lectures evolved with time. {\red As we have seen, March also taught at the ISUP. March had been in 1912 the main initiator of Pearson in France. He translated {\it The Grammar of Science} into French and the book was published in Paris by F\'elix Alcan  (\cite{Pearson1912}). He was probably a very useful source of information on biometry for Darmois.
 Darmois' progressive acquaintance with the tools of the foreign statistical schools (British, American, Scandinavian and so on) convinced him of the necessity of improving the probabilistic level of the lectures. This also perfecly fit Borel's ideas. This rapid evolution led to a reversal of priority between the statistical models and the probabilistic models.  Darmois published the first French textbook on mathematical statistics entitled {\it Statistique Math\'ematique}. The book was published in 1928 by Doin with a preface by Huber. It  shows a strong probabilistic orientation. {\red Darmois eventually obtained a tenured position at Paris University only in 1933. In 1934, he published a a work of popularization about statistics \cite{Darmois1934}.}

A thorough study of the ISUP syllabus is presented in \cite{Roy1937}, \cite{Pressat1987} or \cite{Morrisson1987}. }The comparison between the 1924 and the 1938 syllabi is instructive and gives a clear idea of the state of the statistical scene in France, after 20 years of effort. Now, a large amount of statistics was clearly presented as an emanation of probability theory. The mathematical level subsequently increased and the probabilistic tools introduced during the 1920s and the 1930s wera extensively used. Probability distributions and limit theorems became omnipresent. Probability theory and its {\it main theorems} were introduced at the very beginning, as well as limit theorems and their tools, such as the {\it \v Cebyshev's method for the law of large numbers}. Moreover, a strong presence of analysis of correlation and regression shows that the mentioned notions had become basic knowledge. Furthermore, the number of students of the ISUP had slowly increased: for the academic year 1937-1938, there were 15 students.


\subsection{The Institut Henri Poincar\'e}

As we have seen, the ISUP was a place for teaching, not specifically devoted to mathematics. But in Borel's mind, probability and  mathematical physics also needed a place to develop. That is  why, thanks to the Rockefeller Foundation and Baron E. de Rothschild, in 1928, Borel managed to found an Institute dedicated to probability and mathematical physics, the \textit{Institut Henri Poincar\'e} (Henri Poincar\'e Institute), in short the IHP. The purpose assigned to the Institute was in the first place to facilitate exchanges and meetings between specialists of these domains. 

After the Great War, the Rockefeller Foundation offered credits - both for the victors and the (former) enemy - to develop scientific institutions in Europe through the International Education Board, its organization devoted to this kind of funding. This action not only helped developing new institutions but also saved scientific institutions threatened by the consequences of the war. The foundation paid in particular special attention to the countries where there was a strong cultural academic tradition such as Germany where it originated G\"ottingen Institute for Physics. It also contributed to the development of the Institute for Theoretical Physics in Copenhagen (see \cite{Colasse2002}).  Siegmund-Schultze, in his detailed study of the Rockfeller foundation  \cite{Siegmund-Schultze2001}, gives  details about the creation of the IHP. It is interesting to observe that Borel had direct personal  discussions on the subject with Augustus Trowbridge the head of the International Education Board in Paris.

In his inaugural speech, on 17 November 1928, Borel presented a slightly nationalistic picture of the French origins of probability theory. 

\begin{quote}
Probability calculus and mathematical physics are two sciences whose origin is French to a large extent. Should we mention for probability calculus : Fermat, Pascal, d'Alembert, Buffon, Laplace, Cournot, Joseph Bertrand, Henri Poincar\'e, and for mathematical physics  : d'Alembert, Poisson, Fourier, Amp\`ere, Cauchy and again Henri Poincar\'e? It is only in the second half of the 19th Century that progress in science led us to understand the tight links that exist between both sciences which look distinct at first sight, and to realize that the properties of matter and energy studied by mathematical physics are subject to probability and statistical laws\footnote{Le calcul des probabilit\'es et la physique math\'ematique sont deux sciences dont l'origine est pour une grande partie fran\c caise. Faut-il nommer pour le calcul des probabilit\'es : Fermat, Pascal, d'Alembert, Buffon, Laplace, Cournot, Joseph Bertrand, Henri Poincar\'e, et pour la physique math\'ematique : d'Alembert, Poisson, Fourier, Amp\`ere, Cauchy et encore Henri Poincar\'e? C'est seulement dans la seconde moiti\'e du XIX\`eme si\`ecle que le progr\`es de la science amena  \`a comprendre les liens \'etroits qui existent entre les deux sciences au premier abord distinctes et \`a se rendre compte que les propri\'et\'es de la mati\`ere et de l'\'energie qu'\'etudie la physique math\'ematique, sont soumises \`a des lois de probabilit\'e, \`a des lois statistiques.}.
\end{quote}

In the same speech, Borel described the purpose of the Institute, naming scientists who were supposed to come for talks or lectures at the IHP. He mentioned in particular Einstein and Volterra.

\begin{quote}

Other first-rank scientists have let us hope they will participate. Thanks to them, the Institut Henri Poincar\'e  will really be international, not only by virtue of the students who will attend its lectures, but also because of the professors who will give talks and lectures there. It will not only contribute to the progress of science, but also will act as a bridge between peoples by letting scientists from every country collaborate, meet  and better understand each other\footnote{D'autres savants de premier ordre ont laiss\'e esp\'erer leur concours. Gr\^ace \`a eux, l'Institut Henri Poincar\'e sera c\'eritablment international, non seulement par les \'el\`eves qui suivront son enseignement mais aussi par les professeurs qui y donneront des cours et conf\'erences. Il ne contribuera pas seulement au prog\`es de la science, mais au rapprochement des peuples en permettant \`a des savants de tous les pays de collaborer, de se connaitre eet de mieux se comprendre.}. 
\end{quote}

The opening of the IHP allowed Borel to create some new courses in Paris. Borel's lectures on probability (entitled \textit{Probability theory and its application}) were given at the IHP. Other  lectures on the subject were given by Maurice Fr\'echet whom Borel had managed to call from Strasbourg as a professor {\it without chair} (`sans chaire') at the Sorbonne, with the  obvious objective of  being in charge of probabilities in the new Institute. The lectures Fr\'echet gave were entitled \textit{The law of large numbers} and  \textit{The theory of integral equations}. As was said in Borel's inaugural speech, foreign lecturers were invited every year as well to give talks about their field of interest and investigation.

 The main subject taught at the Institute was however physics, as well as several topics of  {\it applied mathematics}. Darmois was asked by Borel to give lectures at the {\it Institut} on statistics from its very opening, in 1928, when he was still a professor in Nancy. An announcement from academic year 1928-29 asserts that {\it M.G.Darmois, professor at the University of Nancy, will give a series of four lectures on the following topic : Statistical laws, correlation and covariance with applications to heredity, to social and economic sciences.} The lectures took place in April 1929. During the first year, there was no mention in the syllabus of other lectures related to statistics. Until he obtained an academic position in Paris, Darmois would come every year to give some talks and lectures at the IHP and he would remain the only speaker on statistics. The list of invited foreign speakers clearly emphasizes that scientists dealing with physical questions were mainly invited. During the first two years, Fermi, Einstein, Dirac or also Sommerfeld were invited, to mention only these famous names. Among the first specialists in probability, we see Hostinsk\'y, P\'olya, and also Paul L\'evy. The part of probability progressively grew among the subjects  of the conferences - and later statistics was also represented. 
 
 An important feature of the IHP was that Borel used the new institution to help his foreign colleagues in the grip of the European interwar convulsions. The administrative archives of the Institute contain several interesting documents in that respect. For example a letter from the minister for public instruction to Borel (via  Dean Maurain) in 1932, mentioning that Volterra had been { \it disbarred from university staff for political reasons and left for Paris}. After Hitler's nomination as Chancellor in Germany in 1933, an increasing number of German refugees were invited to give talks at the IHP\footnote{Among them, Gumbel settled down in France - see \cite{Hertz1997}.}. The ones Borel never succeeded to invite were the Soviet scientists confronted with closed borders after Stalin's power consolidated. Whereas Moscow was in the 1930s the indisputable world center of the mathematics of randomness, it is rather amazing not to see a single Soviet name among the speakers. The presence on the list of mathematician Vladimir Alexandrovitch Kostitsyn was only due to the fact that he managed to come to Paris in 1928 under the pretext of medical emergency and never returned to the USSR.

As we have already mentioned, apart from Darmois, no statistician was  invited during the first three years of existence of the IHP. Only in 1931 did a first foreign statistician give a talk. Subsequently in the 1930s most statistical trends were represented in these conferences : the Scandinavian (Guldberg, Steffensen), the British (Neyman, Fisher), the German (Gumbel) in particular. This was the result of a kind of deliberate strategy for importing foreign statistical technique to France. We shall now deal with this question.

\section{Technology transfer}
We have already evoked in the previous section how France lagged behind on statistical methods when compared to other countries. We shall try to emphasize the role played by Borel to improve this situation - not the least via the creation of the IHP.

{\red It seems there was a convergence of view between Borel, Darmois and Fr\'echet concerning the statisticians whose works had to be presented at the IHP. In 1933 a volume of Borel's great {\it Treaty of calculus of probability and its application} was published by Risser and Traynard (\cite{RisserTraynard1933}) - the only one whose title included the words {\it Statistique Math\'ematique}. Most of the material of the conferences on statistics at the IHP during those years were included in the book. The American mathematician Edward Lewis Dodd, whom Aldrich claim to be one of the most influent American statisticians of the time  (\cite{Aldrich2007}), wrote a review for the American Mathematical Society where he underlined the up-to-date treatment of mathematical statistics by Risser and Traynard (\cite{Dodd1934}).}

\subsection{The {\it Annales}}

To facilitate the diffusion of the written text of the conferences, the  {\it Annales de l'IHP : Recueil de Conf\'erences et M\'emoires de  Calcul des Probabilit\'es et Physique Th\'eorique} (Annals of IHP : A collection of Conferences and Memoirs on Calculus of Probability and Theoretical Physics) were founded in 1930. Contrary to the current publication bearing the tittle {\it Annales de l'IHP}, the {\it Annales} in the 1930s were not really conceived as a journal but as a publication devoted to recording the lectures held at the Institute. For many years, there were no spontaneous submissions of papers. Our efforts to find announcements of this publication abroad has produced only slight results, which nonetheless reveal at least a limited diffusion. The only true advertisement we noted was in the Bulletin of the AMS 
in 1931 (\cite{AMS_IHP_1930}). It reads 

\begin{quote} A new journal entitled Annales de l'Institut Henri Poincar\'e, Recueil de Conf\'erences et M\'emoires de Calcul des Probabilit\'es et Physique Th\'eorique has been founded by the Institut Henri Poincar\'e for the publications of lectures delivered by invitation at that institute. The first issue contains articles by C.G. Darwin, A. Einstein, and E. Fermi.
\end{quote}

Until 1964, there was only one series of the {\it Annales},  mixing papers on physics or mathematical analysis with others on probability or statistics. In fact, a majority of papers were devoted to physics. The list of publications on probability and statistics reveals 17 articles before World War 2.

With the exception of Paul L\'evy (maybe because his talks were organized by oral agreement), we have been able to locate the moment when the lectures were given. There was generally  an average two years' delay between the lecture and its publication, and this probably resulted in a final form slightly more fully worked out than the oral lecture. As the lectures were really an opportunity for the speaker to present research in progress, the later publication was enriched with more recent results and discussion. Through the politically-correct tone of the academic prose, one may sometimes feel several scientific polemics of the time such as the fiery discussions around the foundations of probability theory (von Mises was invited in 1931, Cantelli in 1932, de Finetti in 1935, Reichenbach in 1937). Or one can also follow the evolution of a new subject such as the theory of Markov chains. Hostinsk\'y's 1937 paper, made after his lectures the same year, testifies to the predominance of the Soviet school of probability with Kolmogorov as a beacon, and also refers to Doeblin's work (see \cite{MazliakDoeblin2007}). 

\centerline{\includegraphics[width=15cm]{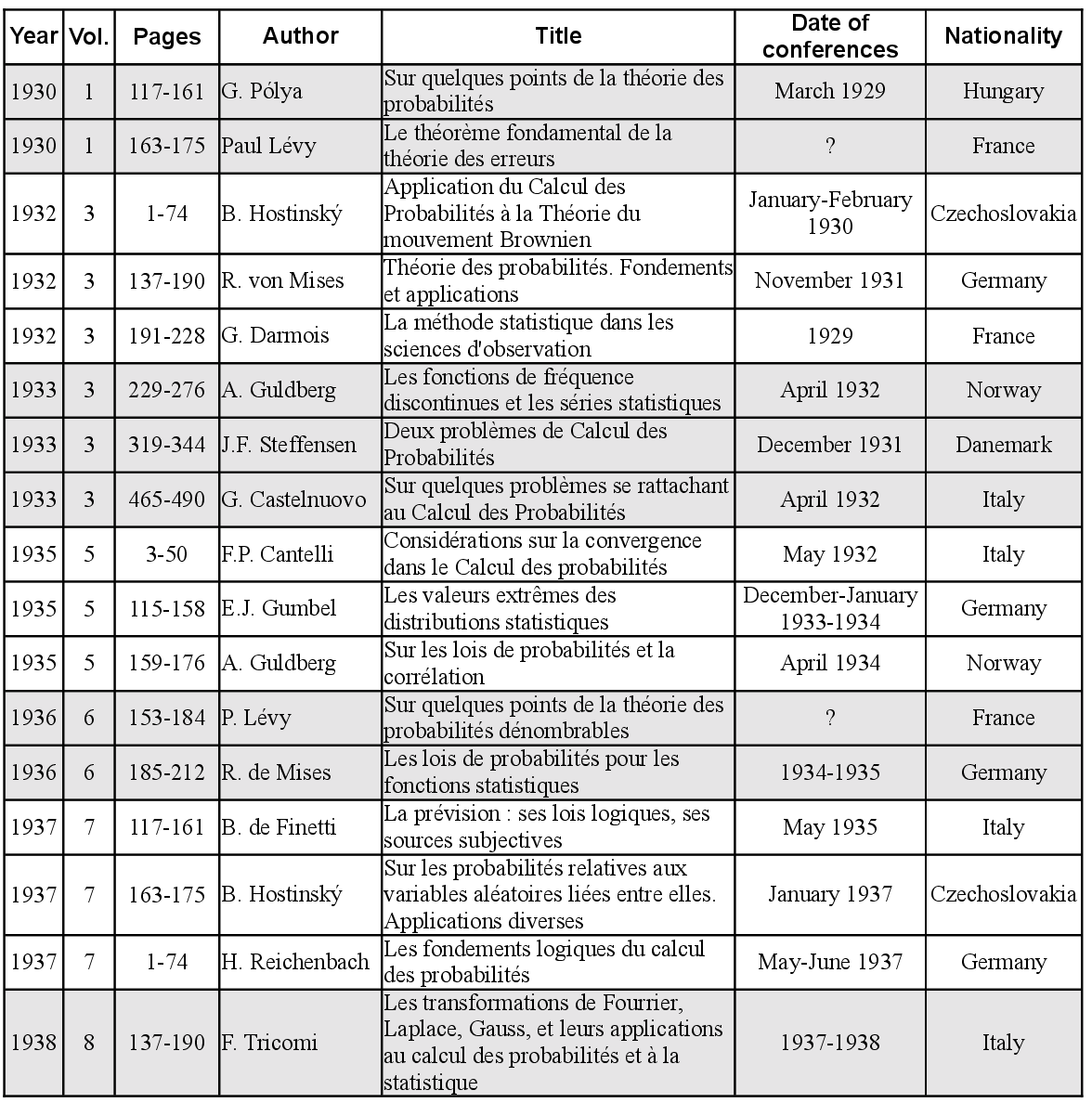}} 
\centerline{Papers about probability and statistics} \centerline{in {\it Annales de l'Institut Henri Poincar\'e} during the 1930s}

\vskip 3mm

 All this draws a picture of the IHP as a place with lively exchanges in the study of mathematics of randomness during the 1930s. The collection of the 17 papers of the {\it Annales} reflects the efforts of the time in that direction - with the noticeable and already mentioned absence of Soviet authors. A certain prevalence of theory over application can be observed - which results in an orientation which favors probability over statistics. This is not surprising,as we said in the previous section that in Fr\'echet and Darmois' vision of the mathematics of randomness, probability theory was the basis on which to build statistics. More significantly, we observe that important foreign traditions are represented in this list:  Scandinavian (Guldberg, Steffensen), British (Neyman), German (Gumbel) ; as we shall see below,  Darmois' main purpose in his 1929 lectures (the first lectures in statistics at the IHP) was to present general methods obtained by Pearson's biometric school in Britain. In fact, a careful observation of the statistical methods exposed at the Institute reveals significant choices. We have mentioned that Borel had known about Pearson's methods early in his scientific life, maybe through Volterra who had always advocated them (see \cite{Volterra1906} and  \cite{Volterra1939}). But we have seen that Borel had also expressed deep criticisms towards  other aspects of the British studies of randomness, {\red when he commented on Keynes' treatise.}

Did Borel's distrust for British {\it different minds} exert an influence on the choices? R.A.Fisher, for example, who became Pearson's successor as the central figure of British statistics during the 1930s, was invited only in 1938 and his lecture was not published. Nor does it seem that his methods were presented at the Institute by other scientists\footnote{It is true that Darmois quotes Fisher twice in his paper (pp. 209 and 223) but it is only in passing, as a source of original applied situations.{\red We note further that Edward L. Dodd, in his laudatory review about Risser and Traynard's book  (\cite{RisserTraynard1933}) did not seem to know Fisher's work either. It is only from 1935 on that Darmois grasped the epoch-making aspect of Fisher's  techniques about estimation. Darmois' lectures  at ISUP in 1937, edited by Pierre Thionet (\cite{Darmois1948}) show the high technical level attained in the exposition of notions such as sufficient statistics and hypothesis testing.} Another explanation for Fisher's absence was that Borel, Darmois and Fr\'echet probably observed the intellectual wars raging in Great Britain between the supporters of and opponents to Bayesian inference and it does not seem that they wished the IHP to become a battlefield for these quarrels. On the complicated story of British statistics of that time (and especially on Keynes' intervention in the debate), see \cite{Howie2004} and \cite{Aldrich2008} and the references included.}.

As an illustration of the technology transfer in statistics sought by the IHP team, we shall now concentrate on the first two papers published in the {\it Annales} devoted to statistical topics. These two papers, written by Darmois and Alf Guldberg, were the occasion for presenting results of Pearson's school, and of the Scandinavian tendency.

\subsection{Darmois' paper}

Darmois took time before finalizing the manuscript of the 1929 lectures mentioned above. The published version \cite{Darmois1932}, which appeared in 1932, indicates that it was received on 25 July 1931. 

Rather than a paper discussing a precise theoretical point, Darmois presented a survey about series of observations, and he presented a collection of methods to detect what he called `persistences' (in French {\it permanences}), that is to say laws, deterministic or not, which provide a theoretical framework for the model under study. The program is clearly set at the beginning of the paper.

\begin{quote}

The aim of these lectures is to emphasize several fundamental ideas of statistics and the character of the special order they can bring into some sets of experimental results.

A research project usually stems from the desire to verify,  deepen and link a group of ideas, either because of its practical importance, or from its unquestionable presence. 

As observation advances in the field under study, the accumulated  results  are thus placed either within the framework of a prior theory, sometimes vague, or in a purely pragmatic order, generally with some logical threads running through it. 

This provisional organization, as a convenient description of the observed material, allows the mind to handle large sets more easily.

After that, one looks for what we may generally call persistences, which is to say relations commonly observed which constitute empirical laws for the observed phenomena. (\dots  Once these empirical laws are obtained,  the subsequent process of the mind aims at understanding them, explaining them, relating them. 

In parallel with the descriptive order, the explanatory or logical orders, scientific theories and laws develop. This development, which is constructed from several basic notions, must then recover the results of the experiment\footnote{Le but de ces le\c cons est de mettre en \'evidence quelques id\'ees essentielles de la statistique et le caract\`ere de cet ordre sp\'ecial qu'elle peut mettre dans certains ensembles de r\'esultats exp\'erimentaux.
Un ensemble de recherche se pr\'esente g\'en\'eralement, soit comme suite \`a un d\'esir de v\'erifier, d'approfondir et lier un groupe d'id\'ees, soit comme impos\'e par son importance pratique, par son indiscutable pr\'esence.
A mesure que l'observation progresse dans le domaine qu'elle a choisi, les r\'esultats qu'elle accumule sont alors plac\'es, soit dans les cadres parfois assez vagues d'une th\'eorie pr\'ealable, soit dans un ordre de pure commodit\'e, g\'en\'eralement sillonn\'e de quelques filaments logiques.
Cette organisation provisoire, description commode du mat\'eriel observ\'e, permet \`a l'esprit de dominer avec plus d'aisance des ensembles \'etendus .
On cherche ensuite, ce que d'une mani\`ere g\'en\'erale, nous appellerons des permanences, certaines relations qui demeurent courante, et qui constituent des lois empiriques des ph\'enom\`enes observ\'es.(\dots )
Ces lois empiriques une fois obtenues, les d\'emarches suivantes de l'esprit visent \`a les comprendre, les expliquer, les relier ntre elles.
Parall\`element \`a l'ordre descriptif se d\'eveloppent l'ordre explicatif, logiques, les th\'eories et les lois scientifiques.
Ce d\'eveloppement, qui se construit \`a partir de certaines notions \`a la base, doit alors retrouver les r\'esultats de l'exp\'erience}.
\end{quote}

The paper can be roughly divided into four parts (which may correspond to the four lectures Darmois presented),  separated by a three-star symbol.

The first section is devoted to a general reflection about the notion of statistical permanence. It is deprived of any mathematical apparatus and essentially seeks to illustrate the notion through classical examples such as Mendel's  hybridization, or radioactivity. Darmois' examples emphasize the role of the urn model, where alternatives are represented by balls of different colors. Darmois insists on the fact that the urn model can be very useful as a tool even if it may not provide an explanatory model. For the example of radioactivity, he explains : {\it Everything happens as if there were the same probability of disintegration for each atom\dots In each atom, there is an urn from which the drawings are made} though, naturally, {\it the explanation seems quite intricate} - an intricacy that Darmois in fact does not absolutely reject due to the {\it new mechanics}  (wave mechanics) for which he refers to Haas' 1928 book \cite{Haas1928}, translated into French in 1930 and prefaced by Borel. But this is not the point in Darmois' statistical presentation and, as he stipulates in his previous example of Mendelian hybridization, to represent the statistical persistences, one may work with the urn considered {\it as a fact}.

The second part of the paper is devoted to the presentation of several urn models, with increasing sophistication to render them more adapted to describe particular situations of {\it `one parameter persistences' (permanence \`a une variable)}. Darmois illustrates the elementary case of how to obtain a confidence interval by means of the de Moivre theorem for a simple Bernoulli (independent and identically distributed) sample. Darmois emphasizes Pearson's observation that the aim is {\it to obtain an image of the statistical persistences from the theory of probability}. This implies that it cannot be sufficient to provide limit theorems to deal with a statistical problem. Only a clever use of these theorems can provide a model for the fluctuations or the stability. This insistence on the necessity of obtaining  a good-quality approximation shows that what nowadays seems an obviously required part of the statistical treatment was at the time not so well realized. Darmois presents successively Poisson, Lexis, Borel, and Polya's urn schemes. Interestingly, he proposes concrete data which, he says,  correctly fit these models. The examples chosen for Lexis and Polya's models are rather classically related to mortality or infection data, the kind of figures which invaded statistical studies in continental Europe as well as in Great Britain. However, when completing his list of models with considerations about the Gaussian distribution, Darmois refers the reader to recent economic works by Gibrat (\cite{Gibrat1930})\footnote{{\red About Gibrat, see \cite{Divisia1932} and \cite{Armatte1995}}}. 
Desrosi\`eres (\cite{Desrosieres2000}, pp.200-203) has already observed that before 1930, the use of statistical models by economists had not been frequent. Darmois, as well as sociologists Halbwachs and Simiand, were innovators through their contact  with economists during the 1920s (such as Divisia\footnote{On Divisia, see in particular \cite{Roy1965}
\cite{Lemesle1998}.} 
(1889-1964) - the same that introduced the word {\it permanence}, or Gibrat) with strong mathematical backgrounds who were sensitive to statistics. At the beginning of the 1930s, all these scholars met at the newly founded {\it Centre polytechnicien d'\'etudes \'economiques} devoted to the study of the Great Depression.

The third part of the paper is devoted to the general question of correlation. After a brief presentation of stochastic dependence and of elementary theory of regression (defining the correlation coefficient and linear regression), Darmois gives general results obtained by Galton and Pearson's biometric school. The main examples deal with hereditary phenomena where the question is to bring precision about the influence of ancestors on their descendants (the presence of a character, sex ratio at birth for examples). As Darmois writes: { \it The aim of the biometric school is to obtain connected distributions; in other words, given the value of a character in the father and the grandfather, to deduce the distribution of the same character in the sons with this very ancestry\footnote{Le but de l'\'ecole biom\'etrique anglaise est en effet d'obtenir des r\'epartitions li\'ees, autrement dit, connaissant la valeur du caract\`ere pour le p\`ere et pour le grand-p\`ere, d'en d\'eduire la loi de r\'epartition du m\^eme caract\`ere pour les fils ayant cette ascendance fix\'ee.}}(\cite{Darmois1932}, p.213).

The last and fourth part of Darmois' paper introduces time as a parameter and studies the specific kind of models which can be developed to face this situation. As Darmois explains, time  was naturally present for example in the case of astronomic observations, but with {\it the particular chance that the motion of the earth and the planets around the sun produce simple and rather pure regularities due to the enormity of the mass of the sun.  ( \dots ) But those problems presented by social and economic sciences that statistical methods may seek to clarify and to solve, are very different\footnote{Avec cette chance particuli\`ere que le mouvement de la terre et des plan\`etes autour du soleil fait appara\^\i tre des r\'egularit\'es simples et assez pures, dues \`a la masse \'enorme du soleil.
(\cite{Darmois1932}, p.218 )}}.
The evolution of a quantity is therefore presented in this fourth part as a model integrating random perturbations. Darmois presents several models of that kind, but mostly limits his presentation to description and commentaries, not even summarizing the corresponding technicalities. This limitation to a qualitative exposition  is probably interpretable as  another sign of the low level of statistical knowledge in France at the time, making it necessary to provide a general survey of existing situations before considering the mathematical treatment. Moreover, the four papers referred to in this part were recent publications from the Anglo-Saxon world, a trace of  Darmois' efforts during the 1920s to get acquainted with British and American works in that direction. 

The first model described by Darmois was considered by Udny Yule  in his 1927 paper \cite{UdnyYule1927} in relation to the observations of sunspots. It concerned the presence of random disturbances of an oscillatory system which Yule illustrated by the image of a pendulum {\it left to itself, and unfortunately boys\footnote{not girls\dots The detail may be of importance!} got into the room and started pelting the pendulum 
 with peas.} 
 Yule proposed two methods to investigate the hidden periodicities of the phenomenon, but as we mentioned, Darmois did not  go into such details.
 Darmois found another kind of model in Hotelling (\cite{Hotelling1927}). Hotelling considered a differential equation whose coefficients contains a random perturbation. 
 Hotelling explained in his paper how it may be more appropriate to work with the coefficients of the equation rather than with the solution of the equation when a statistical verification of the model is needed. However, also here, Darmois limits himself to mere description. 
Other models are borrowed  from Ronald Fisher and Egon S. Pearson's recent papers.
In conclusion, Darmois' paper appears as a vast catalogue of situations opening the door to future mathematical studies, with a noticeable inclination towards economic phenomena. 


\subsection{Guldberg's paper}

Let us now turn towards the second paper, written by Norwegian statistician Alf Guldberg (1866-1936). 

Few documents are easily accessible on Guldberg. {\red About his youth as a Lie student traveling in Europe, one may consult \cite{Stubhaug2002} chapter 24. Back in Norway in 1903, he obtained a position at Christiania\footnote{The name of Oslo before 1925.} University and soon began to teach actuarial statistics. Guldberg was often mentioned  in the \v Cuprow-Bortkiewicz correspondance (\cite{Sheynin2005}) as an influential member of the Scandinavian mathematical scene}.

At the beginning of the 20th Century, Scandinavia was a first-rank center for approximation techniques, not behind the British school (see \cite{Schweder1999}). {\red Those techniques were inherited from the works of astronomers of the second half of  19th century, such as Thiele (1838-1910), Gyld\'en (1841-1896) or Charlier (1862-1934)\footnote {See \cite{Mazliak2009} for the filiation between the statistical approach of the continuous fraction development proposed by Gyld\'en and Borel's theory of denumerable  probabilities, and also \cite{Lauritzen2002} and \cite{Holmberg1997} about astronomers in Scandinavia.}\dots} In a draft of an introduction for Steffensen, when he came to Paris in 1932 to give a talk at the IHP,\footnote{We found this draft by chance among the incredible mountain of documents of Fr\'echet's archives in the Paris Academy of Sciences.},  Fr\'echet declared

\begin{quote}
I have the great pleasure to introduce to you M.Steffensen, professor at the University of Copenhagen. It will be for us a new experience which, I am convinced,  will be very valuable. Up to now, the lectures which have dealt with probability calculus have only considered it on the theoretical plane or through its application to physics. M.Steffensen, to whom the whole theoretical part of probability calculus is familiar, is now well known for the way he wonderfully  had allowed the actuarial sciences to profit from mathematical advances. In a book about interpolation  he systematically introduced a demand too often left aside : never to introduce an approximation without trying to specify the error. It is a relatively moderate advance in theoretical mathematics but it did not penetrate deeply into applications. Thanks to M.Steffensen no one can any longer ignore first that this estimation of error is necessary, and second, that thanks to him it is feasible.\footnote{J'ai un tr\`es grand plaisir \`a vous pr\'esenter M Steffensen professeur \`a l'Universit\'e de Copenhague. Ce sera pour nous une nouvelle exp\'erience et une exp\'erience qui j'en suis sur sera tr\`es heureuse. Les conf\'erences qui vous ont parl\'e jusqu'ici du calcul des probabilit\'es ne le montrent que du point de vue th\'eorique ou du point de vue des applications aux Sciences Physiques.
M. Steffensen \`a qui toute la partie th\'eorique du calcul des probabilit\'es est famili\`ere, s'est fait conna\^\i tre par la fa\c con magistrale dont il a su faire profiter les sciences actuariels des
progr\`es math\'ematiques. C?est ainsi que dans un livre sur l'interpolation il a introduit syst\'ematiquement une exigence trop souvent laiss\'ee de c\^ot\'e : celle qui consiste \`a n'introduire aucune forme approch\'ee sans chercher \`a en chiffrer l'erreur. C'est l\`a un progr\`es relativement mod\'er\'e accompli dans la math\'ematique th\'eorique mais qui n'avait pas p\'en\'etr\'e profond\'ement ses applications. Gr\^ace \`a M. Steffensen on ne peut plus ignorer d'abord que cette estimation de l'erreur est n\'ecessaire, ensuite que gr\^ace \`a lui elle est possible.}
 \end{quote}

Steffensen's book  referred to by Fr\'echet, \textit{Interpolationsl\ae re} was published in Danish in 1925, and  an American translation  (\cite{Steffensen1927}) rapidly became available, appearing in 1927. In the Bulletin of the AMS \cite{SteffensenAMS1928}, the notice explained what Steffensen's book contains.

\begin{quote}
Professor Steffensen's book is the outgrowth of the lectures which the
author has given to actuarial students at the University of Copenhagen
and is, with a few additions and simplifications, a translation of the Danish
edition published in 1925. The book is intended as a text for students in
American colleges and requires as mathematical equipment only an elementary
knowledge of the differential and integral calculus. In a few places
where the gamma function has been used the paragraphs have been printed
in smaller type and may be omitted without breaking the continuity of the
text.
The topics covered are (1) the general theory of interpolation and extrapolation
including the standard formulas usually associated with the names
of Newton, Gauss, Stirling, Bessel, and others; (2) numerical differentiation;
(3) numerical integration; (4) numerical solution of differential
equations.
Professor Steffensen's treatment is more rigorous than is usual in books
on interpolation. This is important not merely from the point of view of
the pure mathematician but also because of the increased number of formulas
with workable remainder terms. It should not be supposed, however,
that this adds to the difficulty of reading the text. The style is clear and,
after the meaning of the symbols has been mastered, the book should
prove very valuable to the increasing number of Americans who require
some knowledge of this field of mathematics. The formulas and methods
are illustrated by simple numerical examples, but the value of the book
for class room use would be increased if it contained some problems to be
solved by the student. 
\end{quote}

The mathematical study of interpolation was Steffensen's known speciality, as \v Cuprov in Prague  wrote to Bortkiewizc  \footnote{See \cite{Sheynin2005}, letter dated 2 April 1925.}:
{\it Steffensen's new textbook is once again written in Danish. It contains considerations on interpolation, which is Steffensen's speciality.
}
The American version remain the only translation of the book. In France, it seems to have been unknown for quite a while.

More important for Guldberg's lectures in Paris was another aspect of Scandinavian statistics. Thiele had defined in the 1880s a mathematical object, the \textit{halfinvariants},  later rediscovered by Fisher under the name {\it cumulative moment functions} and finally \textit{cumulants}, as they are called today. {\red We shall comment on this notion later}.

As we mentioned above when we described Darmois' paper, in the 1930s, the British statistical school was predominant in Western Europe and Guldberg's presence in Paris may have been a reaction to advertise the methods of the Scandinavian school.  Two points could corroborate this hypothesis. Fisher's paper - where he introduced cumulants - was published in 1929 \cite{Fisher1929}, and there was no reference to Thiele in it.  Fisher was later  informed about Thiele's priority. He nevertheless did not entirely accept the fact and constantly underrated Thiele's treatment of halfinvariants\footnote{See details in Lauritzen \cite{Lauritzen2002}, or Hald \cite{Hald1998} (pp.344-349).}.
In Guldberg's paper (\cite{Guldberg1932}), there is clearly a desire to display Scandinavian works, and Thiele's in particular. More generally, Guldberg seems to have engaged in frenetic activity on the international scene to obtain recognition of Nordic  mathematics. In particular, he was the main artisan of the Oslo International Congress of Mathematician in July 1936, according to Carl St\o rmer, president of the congress, who replaced Guldberg  at the last moment due to the latter's sudden death in February.

\begin{quote}

It was during the last international congress in Z\" urich that the Norwegian mathematicians took the decision to propose that the next congress could be held in Oslo\footnote{In \cite{Zurich1932}, the report of the final session mentions that indeed Guldberg made the suggestion which was acclaimed by the Congress.}. It is our colleague, the regretted Alf Guldberg, who took this initiative, and who since then had constantly worked to guarantee the success of the congress, and above all to place it on a reliable economic basis. We had hoped that Guldberg could open the Congress as president. But, as every member of the congress knows well, a premature death snatched him away from his work. Nobody but we, who must finish his work, can better realize what a loss it was for the Congress to have been deprived of our eminent friend, who was so efficient in attracting friends and knew the art of association with his  fellows ({\it sic})\footnote{Ce fut au dernier congr\`es international de Zurich que les math\'ematiciens norv\'egiens r\'esolurent d'emettre le projet que le prochain congr\`es se r\'eunirait \`a Oslo. C'est notre coll\`egue, le regrett\'e Alf Guldberg, qui en orit l'initiaive, et qui depuis s'est pr\'eocupp\'e inlassablement \`a assurer la succc\`es du congr\`es, et surtout \`a l'\'etablir sur une base \'economique solide. Nous avions esp\'er\'e que Guldberg pourrait ouvrir le Congr\`es comme pr\'esident. Mais, comme le savent tous les membres du congr\`es , une mort pr\'ematur\'ee vint l'aracher \`a son oeuvre. Personne mieux que nous qui avons d\^u achever son travail, ne peut savoir quelle perte ce fut pour le congr\`es que d'\^etre contraint \`a se passer de notre \'eminent ami, qui savait si bien gagner des amis et connaissait l'art de fr\'equenter ses semblables.(\cite{Oslo1936})}.

\end{quote}

{\red
Let us briefly introduce Thiele's halfinvariants. Thiele first presented the halfinvariants in his 1889 book  \textit{Almindelig Iagttagelsesl\ae re: Sandsynlighedsregning og mindste Kvadraters Methode} (The general theory of observation: Calculus of probability and the method of least squares),  translated and edited in \cite{Lauritzen2002}. 

If we denote by $\sigma_{n}$ the $n$-th moment of a random variable, the halfinvariants $\mu_{n}$ were first defined by the following recursion formula

$$\sigma_{k+1}=\sum_{i=0}^{r} 
\left( \begin{array}{c}
k\\i
\end{array} \right)
 \sigma_{k-i}\mu_{i+1}$$

Moreover, Thiele explained in his book that the first four halfinvariants were the most important ones. The first was equal to the mean, the second to the variance, the third measured the {\it skewness} and the fourth the flatness of a distribution. For Hald (\cite{Hald1998},  p.209) there are basically two important questions.

\begin{quote}
\begin{enumerate}
\item How did Thiele find the recursion formula defining the cumulants?

\item Why did he prefer the cumulants over the central moment?
\end{enumerate}
\end{quote}

In his book of 1889, Thiele explained that the halfinvariants are often small.

\begin{quote}
We shall later become acquainted with certain favourable properties that distinguish the halfinvariants from the reduced sums of powers; here we only mention that the halfinvariants most often are numerically smaller than the others. Whereas the reduced sums of powers of even degree are always positive, all halfinvariants except $\mu_{2}$ may be negative just as well as positive, and therefore they will often have values that are not far from zero. (\cite{Lauritzen2002}, pp. 84-85)
\end{quote} 

Only ten years later, in paper \cite{Thiele1899} did he realize that the halfinvariants could be directly expressed as coefficients of the Taylor expansion of the Laplace transform's logarithm. 
$$
{\rm exp}(\sum_{i=1}^{\infty} \frac{\mu_{i}t^{i}}{i!} )
=
\int e^{tx}f(x)dx .
$$

In \cite{Thiele1899}, Thiele proved that the halfinvariants for the Gaussian distribution posess an important property: only the  first two cumulants are non equal to zero. This could have been another fundamental reason to use them. }

Guldberg went to the IHP in April 1932 for the first time. He returned later to the IHP in April 1934. His two talks were published in the {\it Annales } as \cite{Guldberg1932} and \cite{Guldberg1934}. Both papers mainly deal with the same statistical technique, namely the application of Thiele's halfinvariants to the identification of an unknown probability distribution: in \cite{Guldberg1932} for a single random variable, in \cite{Guldberg1934} for two. Let us now comment on \cite{Guldberg1932}, which was preceded in 1931 by a Note to the Comptes-Rendus (\cite{Guldberg1931}), presented by Borel.

The article is divided into four parts. In the first one, Guldberg introduces definitions and notations, and applies them to some examples. In particular, he introduces Thiele's halfinvariants, giving both definitions successively provided by Thiele as we mentioned above. 

\begin{quote}

The reason why Thiele introduced the halfinvariants, a notion which seems rather artificial, is the following : if the observations follow a Gaussian distribution, every halfinvariant of order $r>2$ is 
equal to 0. As one of Thiele's former students, my colleague M.Heegard, made me observe, one has therefore a way for examining if a statistical collection of data can be represented by a Gaussian distribution or not.

Halfinvariants also have other advantages. When one linearly transforms the observations 
$$x_{i}=ax_{i}+b$$
that is to say, when one changes the origin of the observations and the unity which measures them, the halfinvariants are changed according to the following formulae 
$$\mu_{1}^{'}=a\mu_{1}+b$$
$$\mu_{r}^{'}=a^{r}\mu_{1}, r>1$$
(\dots ) The mean moments are transformed in a more complicated way.

These relations show the great importance of halfinvariants for the study of a statistical collection of data\footnote{La raison pour laquelle THIELE a introduit les semi-invariants, notion qui semble quelque peu artificielle, est la suivante : si les observations suivent une loi de GAUSS, tous les semi-invariants d'ordre $r>2$ sont nuls. On a donc, comme me l'a d'ailleurs fait remarquer un ancien \'e\`eve de THIELE, mon coll\`egue M. HEEGAARD, un moyen  d'examiner si un ensemble statistiques peut se repr\'esenter par l'interm\'ediaire de la loi de Gauss ou non.
Les semi-invariants ont cependant d'autres avantages. Quand on fait une transformation lin\'eaire des observations
$$x_{i}=ax_{i}+b$$
c'est-\`a-dire, quand on change l'origine des observations et l'unit\'e qui les mesure, les semi-invariants se tranforment d'apr\`es les formules
$$\mu_{1}^{'}=a\mu_{1}+b$$
$$\mu_{r}^{'}=a^{r}\mu_{1}, r>1$$ (\dots )
Les moments moyens se tranforment d'une mani\`ere plus compliqu\'ee. Ces relations montrent la grande importance que pr\'esentent les semi-invariants pour l'\'etude d'un ensemble statistique}.
\end{quote}

In the second part of his article, Guldberg presents how Pearson represented statistical series, and criticizes  the `British method' because  {\it M.Pearson's curves} are reasonably suitable only in the case when one {\it a priori} knows that the data can be represented by one of these curves\footnote{\cite{Guldberg1932}, p.246}. Guldberg therefore opposes what he calls the  {\it continental method} to the methods of the British biometric school. Let us observe that the  continental statisticians that Guldberg named are almost only Scandinavian: Gram, Thiele, Bruns, Charlier, to which only \v Cuprov and Bortkiewizc were added\footnote{Both kept tight contact with the Scandinavian group during the 1920s
as attested by their correspondence (\cite{Sheynin2005}).}. 

In fact, the real purpose of the paper is presented later, when Guldberg write that he {\it will try to study these questions in another way 
}. He first lists the problems , and then tries to answer to them.

\begin{quote}
The theory of functions of frequencies presents four important problems :
\begin{itemize}
\item The numerical computation of the function
\item How to proceed if one wants to substitute a continuous function for the discontinuous one, which takes precisely the same value as the latter for integer values of the variable
 \item The determination of the moments of a given function of frequencies
 \item A statistical series being given, look for a function of frequencies which gives an approximate representation of it and specify, if possible, the necessary and sufficient criteria for a definite function to fulfill the required conditions\footnote{La th\'eorie des fonctions de fr\'equence pose quatre probl\`emes assez importants \`a savoir :
\begin{itemize}
\item Le calcul num\'erique de la fonction
\item La mani\`ere de proc\'eder si l'on vaut substituer \`a la fonction discontinue une fonction continue, qui pour des valeurs enti\`eres de la variable prenne justement les m\^eme valeurs que la premi\`ere.
\item La determination des moments d'une fonction de fr\'equence donn\'ee[...].
\item Une s\'erie statistique \'etant donn\'ee, chercher une fonction de fr\'equence qui en donne une repr\'esentation approch\'ee et \'eablir, s'il est possible, les crit\`eres n\'ec\'essaires et suffisants pour qu'une fonction d\'etermin\'ee remplisse les conditions equises.
\end{itemize} 
 }
 \end{itemize} 
\end{quote}
The first and the second problems correspond to Guldberg's desire to interpolate the (discrete) empirical frequency function by a {\it continuous} (in the sense of {\it regular} : he asks for derivatives to exist) one, in order to use approximation techniques provided by real analysis. 
The meaning of the third point is clear: it is necessary to compute empirical and theoretical halfinvariants in order to compare them. The fourth point is the general method Guldberg wants to present for the statistical inference of a series of data, a method which relies on halfinvariants. After this general presentation, Guldberg looks at these four questions in particular cases for the theoretical distribution.  For each theoretical distribution he considers (Poisson, Binomial, Pascal and hypergeometric), he provides a analytical expression $\alpha (k)$ depending on the halfinvariants and which is constant. The method corresponding to the last point of his program is then to check whether this function on the empirical halfinvariants is constant. however he only illustrates his full program for the Poisson distribution (\cite{Guldberg1932},  p.235).

\section*{Conclusion}

Whereas mathematical statistics based on probability theory was nearly absent from the French academic scene before the Great War, it had become, after World War 2, a topic well established in France. And it was not only the institutions created in the interwar period (ISUP in Paris, other various institutes in the provinces - such as the IEC - Institute for Commercial studies - in Strasbourg or ISFA  (Institute of Financial Science and Insurance) in Lyon - and the IHP for the theoretical aspects\dots) carried on activities which included this large part of the mathematical sciences; other important structures were founded later, such as INSEE or ENSAE, to facilitate the teaching of specialized sides of those techniques.   As we have seen, it was under the influence of Borel that this field of mathematics began to mature in the 1920s. Borel's activity mixed political and social engagements with his mathematical interests,  for example through the creation of  the journal {\it Revue du mois} he had founded with his wife in 1905. It was during this period of intense reflection about the application of mathematics that Borel realized the role that the measure of randomness would be called upon to play in the future development of various sciences. Borel felt he was in possession of new tools to face mathematical questions, in particular the new analytical techniques of measure theory and Lebesgue integration. Also, risk quantification of life accidents was for Borel an essential information source for the organization of public social institutions : this was the core of the radical political program he was very close to.  Borel therefore became one of the first mathematicians to engage in a renewal of probability theory. With the Great War, and his many involvements in the war effort, he began a career as a politician and used his position to become an active instigator of statistics based on probability. He thought that Paris needed an institute devoted to statistics, and convinced Lucien March and Fernand Faure to create such an institution, the ISUP.  Georges Darmois started to teach in this institute. The ISUP was however only a place for lectures and in Borel's mind, probability and statistics also deserved a real research center. This led to the foundation of the Institut Henri Poincar\'e. Darmois and Fr\'echet, who left Strasbourg after Borel asked him to join the IHP team, decided to use the Institute as a hub where it would be possible to organize a kind of technology transfer for the importation of the methods of mathematical statistics that had been developed abroad (Britain, Scandinavia, Italy\dots). Darmois in particular was eager to introduce British mathematical statistics in France, firstly from Pearson's biometric school, and, at the end of the 1930s, from Fisher and his followers.  

{\bl We tried to describe how the aforementioned  transfer to the IHP was organized, from places where reflection on the mathematization of statistics had been pushed much farther than in France. It would be interesting to analyze more systematically the extent of simple reproduction and of transformation of the contents in this transfer along the lines proposed some years ago by Pestre about the reproduction of the experiments in physics. Pestre writes  that {\it even in the case where [the scientist] chooses to reproduce an experiment to take over it and `check' it, the reproduction is generally `improved' and not ` identical' - and this often leads to complex debates about the relationship between these various experiments and the results proposed by them\footnote{M\^eme dans le cas o\`u il choisit de reproduire une exp\'erience pour se la r\'eapproprier et la "v\'erifier", la reproduction est en g\'en\'eral "am\'elior\'ee", et non "\`a l'identique" - ce qui conduit souvent \`a des d\'ebats complexes sur les rapports qu'entretiennent ces diverses exp\'eriences et les r\'esultats qu'elles proposent ( \cite{Pestre2006}, p.126).}. In our case, what should be focused on is the subsequent path of mathematical techniques with Anglo-Saxon or Scandinavian sources, to better understand how they became integrated (or not) into the mathematical statistics taught and `produced' in France.}

Scrutinizing the {\it Annales de l'Institut Henri Poincar\'e} gives a good picture of how probabilistic statistics appeared in France between the two World Wars when German, Scandinavian and British statisticians came to the IHP to hold conferences. Darmois' and Guldberg's papers, which we have chosen as a sample because they were the first papers in statistics, provide a good illustration of the initial situation. A new generation of probabilists and statisticians obtained a doctorate at the IHP under the direction of Fr\'echet and Darmois at the end of the 1930s. Among them, two were more specifically oriented towards statistical topics : Mal\'ecot, whose research in biostatistics and genetics was based on Fisher's techniques, and Dugu\'e, who became after the war the leader of French statistics at Paris University. 

\section*{Acknowledgement} During the writing of the paper, the authors benefited from two short research stays at the CIRM in Marseille and at the Roscoff conference center of University Pierre et Marie Curie (Paris 6).   We are glad to acknowledge these institutions for the interesting opportunities.

\end{document}